\documentclass[11pt,leqno]{article}
\usepackage{amsthm,amsfonts,amssymb,amsmath,oldgerm,cases,upgreek}
\usepackage{epsfig}
\numberwithin{equation}{section}
\usepackage[thinlines]{easybmat}

\def\eps{\varepsilon }

\newcommand\R{\mathbb R}

\newcommand\C{\mathbb C}

\def\eps{\varepsilon}

\newcommand\br{\begin{remark}}
\newcommand\er{\end{remark}}
\newcommand\bp{\begin{pmatrix}}
\newcommand\ep{\end{pmatrix}}
\newcommand\be{\begin{equation}}
\newcommand\ee{\end{equation}}
\newcommand\ba{\begin{equation}\begin{aligned}}
\newcommand\ea{\end{aligned}\end{equation}}

\newcommand{\bap}{\begin{app}}
\newcommand{\eap}{\end{app}}
\newcommand{\begs}{\begin{exams}}
\newcommand{\eegs}{\end{exams}}
\newcommand{\beg}{\begin{example}}
\newcommand{\eeg}{\end{exaplem}}
\newcommand{\bpr}{\begin{proposition}}
\newcommand{\epr}{\end{proposition}}
\newcommand{\bt}{\begin{theorem}}
\newcommand{\et}{\end{theorem}}
\newcommand{\bc}{\begin{corollary}}
\newcommand{\ec}{\end{corollary}}
\newcommand{\bl}{\begin{lemma}}
\newcommand{\el}{\end{lemma}}
\newcommand{\bd}{\begin{definition}}
\newcommand{\ed}{\end{definition}}
\newcommand{\brs}{\begin{remarks}}
\newcommand{\ers}{\end{remarks}}

\newtheorem{theo}{Theorem}[section]

\newtheorem{exams}[theo]{Examples}

\numberwithin{equation}{section}

\newcommand{\const}{\text{\rm constant}}

\newcommand{\sgn}{\text{\rm sgn}}
\newtheorem{theorem}{Theorem}[section]
\newtheorem{proposition}[theorem]{Proposition}
\newtheorem{corollary}[theorem]{Corollary}
\newtheorem{lemma}[theorem]{Lemma}
\newtheorem{definition}[theorem]{Definition}

\newtheorem{example}[theorem]{Example}
\newtheorem{remark}[theorem]{Remark}

\title{
Existence and stability of steady noncharacteristic solutions on a finite interval
of full compressible Navier--Stokes equations
}

\author{\sc \small
Blake Barker\thanks{
Brigham Young University, Provo, UT 84602;
blake@mathematics.byu.edu: Research of B.B. was partially supported under NSF grant no. DMS-140087.
},
Benjamin Melinand\thanks{ CEREMADE, CNRS, Universit\'e Paris-Dauphine, Universit\'e PSL, 75016 PARIS, FRANCE;
melinand@ceremade.dauphine.fr.
}
and
Kevin Zumbrun\thanks{Indiana University, Bloomington, IN 47405;
kzumbrun@indiana.edu: Research of K.Z. was partially supported
under NSF grant no. DMS-1700279.
 }}
\begin{document}

\maketitle


\begin{center}
{\bf Keywords}: Steady solutions, gas dynamics, Evans function.
\end{center}

\begin{abstract}
We treat the 1D shock tube problem, establishing existence of steady solutions of full (nonisentropic) 
polytropic gas dynamics with arbitrary noncharacteristic data.
We present also numerical experiments indicating uniqueness and time-asymptotic stability of such solutions.
At the same time, we give an example of an (artificial) equation of state possessing a convex entropy for which there holds 
nonuniqueness of solutions. This is associated with instability and Hopf bifurcation to time-periodic solutions. 
\end{abstract}

\section{Introduction}\label{intro}
In this paper, continuing investigation in \cite{MelZ} of the isentropic case, 
we study by a combination of analytical and numerical techniques the existence, uniqueness, and stability of steady solutions of the full (nonisentropic) 1D compressible Navier--Stokes equations on a bounded interval, with noncharacteristic inflow-outflow boundary conditions, and more generally of hyperbolic-parabolic systems of conservation laws of similar abstract type.

This corresponds to the 1D version of the ``shock tube'' problem of describing flow in a finite length and width channel,
with prescribed boundary conditions at the left and right ends. Our main interest is in {\it large-amplitude} data.
Small-amplitude 1D existence, uniqueness, and spectral stability are shown
for general symmetrizable systems in \cite{paper2}.

As developed in the viscous shock case \cite{BHZ,BHLynZ2,BHLynZ3,HLyZ1,HLyZ2}, a convenient method to study spectral stability is via numerical Evans function investigations.
A useful necessary condition, also based on Evans function considerations,
is positivity of the stability index, a mod two count of the Morse index of the linearized operator about the wave.
This was trivially evaluable in the isentropic case \cite{MelZ}, but is complicated 
in general. In particular, it does not seem to be analytically evaluable for the nonisentropic case considered here.
We carry out here both full Evans function and stability index calculations at the same time, both using
the numerical code STABLAB \cite{Barker_matlab}.

\subsection{Description of main results}\label{s:results}
Our main analytical result is the global existence of steady solutions of
the full polytropic gas equations \eqref{nsint} (Theorem \ref{maincor}), proved by a Brouwer degree argument using detailed and special ODE estimates,
applied to the ``Cauchy-to-boundary value'' map $\Psi$ defined in Section \ref{s:profile}.
We show, moreover, that global uniqueness of solutions of \eqref{nsint} is roughly equivalent to transversality of
steady profiles as solutions of the ODE connection problem \eqref{prof}-\eqref{fund}. 
This is equivalent to the nonvanishing of the Jacobian $\det (d\Psi)$ of \eqref{fund} (Proposition \ref{uniqueprop}).

Nonvanishing of $\det (d\Psi)$ is also seen to be equivalent to nonvanishing of the stability index (Lemma \ref{ZSlem}). Hence a change in sign implies appearance of both nonuniqueness and instability: 
the usual ``exchange of stability'' scenario familiar from finite-dimensional ODE.
Thus we may study uniqueness in passing, in the course of a larger study of spectral stability.

Augmenting our analytical results for the full polytropic gas equations, we carry out such a study in Section \ref{s:num} by a systematic numerical 
Evans function investigation of the ``feasible set'' $\mathcal{C}_{u_{0},e_{0}}$ of profiles realizable by numerical shootings. 
Our numerical findings are that, on the feasible set $\mathcal{C}_{u_{0},e_{0}}$,
the stability index is uniformly positive, indicating {\it uniqueness of large-amplitude solutions}, 
and that steady solutions exhibit {\it uniform spectral} {\it stability}. 
We note that nonlinear stability can be shown to follow from spectral stability by similar considerations to 
those of \cite[Section 6]{MelZ}; see \cite{paper2}.

On the other hand, we show numerically in Section \ref{s:nonuniqueness} 
that both uniqueness and stability can fail for gas dynamics with an artificial convex equation of state.


\subsection{Discussion and open problems}\label{s:discussion}

The first local existence/uniqueness result for small-amplitude data 
was established in \cite{KK}, in multi-D.  We extend
that result here to large-amplitude data in the 1-D case by a combination of analytical (existence) and
numerical (uniqueness) investigation, establishing (numerically) time-evolutionary stability as well.

Our findings of global existence and uniqueness for the noncharacteristic 
problem parallel those of Lions \cite{Lions_book} in the characteristic case $u=0$ on the boundary,
for which he shows global existence and uniqueness of solutions for arbitrary prescribed average density, 
in 1- and multi-D.
However, they are obtained by quite different techniques, which, moreover, are special to 1D.  Indeed, though perhaps intuitively expectable,
especially given the uniform shock stability results of \cite{HLyZ1,HLyZ2} for the  compressible Navier--Stokes equations, 
our results of large-amplitude existence, uniqueness and stability are obtained by
a combination of exhaustive numerical investigations, and rather delicate degree-theoretic arguments specific to the
equations of 1D polytropic gas dynamics under study.

Our investigations of stability belong, rather, to a newer family of investigations blending numerical and analytical
techniques to study dynamics and bifurcation of shock waves and related solutions of hyperbolic-parabolic 
conservation and balance laws, cf. \cite{HLyZ1,BFZ,BHLZ15,BMZ,SZ,Z1}.

In \cite{SZ,Z1} the case of steady solutions on a half-line was investigated and it was shown that instability of steady solutions can occur, even for the most standard ideal polytropic gas law. 
This suggests that the question of stability at least is not a foregone conclusion for steady solutions on the interval.
Moreover, the nature of instablity found in \cite{SZ,Z1} involved change of sign in the stability index,
which in the present case would signal nonuniqueness as well.
On the other hand, our numerical findings (Section \ref{s:num}) indicate that neither of these phenomena in fact occur for polytropic gas dynamics
on the interval.

This begs the question whether such detailed and special arguments are necessary, or whether there might
instead exist some more straightforward argument for all or part of our results
via general principles, such as, e.g., existence of convex entropy as used in Section \ref{s:existence}.
We give a partial answer to this question in Section \ref{s:ceg}, exhibiting a counterexample
involving an equation of state presented in \cite{BFZ} for which the equations of compressible gas dynamics 
possess a convex entropy, but global stability and uniqueness are violated.
It is seen that the associated transition to instability can involve either steady bifurcation to multiple solutions,
or {\it Hopf bifurcation to time-periodic solutions}.
The latter phenomenon is significant as the first example of Hopf bifurcation for stationary solutions
of compressible gas dynamics, similar to ``galloping'' or ``cellular'' instabilities in detonation
\cite{TZ}.

It is an interesting question whether our existence result extends to general equations of state considered in \cite{BFZ}. Note that we obtain nonuniqueness results for a particular equation of states in Section \ref{s:ceg}.
A further very interesting open problem is the extension of our 
large-amplitude existence results to the true multi-D shock tube problem,
generalizing the small-amplitude existence-uniqueness results of \cite{KK}, 
and the determination of stability of steady multi-D solutions even in the small-amplitude case.



\section{Preliminaries}\label{s:prelims}
\subsection{Equations of motion}\label{s:motion}
The 1D compressible Navier--Stokes equations in Eulerian coordinates are
\ba\label{nsint}
\rho_t+(\rho u)_x &=0\,,\\
(\rho u)_t+ (\rho u^2+p)_x &=\alpha u_{xx}\,, \\
(\rho E)_t+(\rho u E+pu)_x &=\kappa T_{xx}+(\alpha uu_x)_x\,
\ea
where    
\begin{equation*}
E=e+\frac{u^2}{2}, \quad p=\Gamma\rho e,\quad e=c_vT\,,
\end{equation*}
and $\nu = \frac{\kappa}{c_v}$. Here $\Gamma$, $c_v$, $\nu$ and $\alpha$ are fixed positive constants; see \cite{Ba,HLyZ1,HLyZ2}. 

As described in \cite{MelZ} in the isentropic case,
we seek steady solutions on the interval $[0,1]$, 
with noncharacteristic inflow-outflow boundary conditions
\be\label{BC}
(\rho, u,e)(0)=(\rho_0, u_0,e_0),\quad
(u,e)(1)=(u_1,e_1). 
\ee
By changing $\rho$ by $\rho_{0} \rho$,  $u$ by $\frac{1}{\rho_{0}} u$, $t$ by $\rho_{0} t$ and $e$ by $\frac{1}{\rho_{0}^{2}} e$ (notice that we can not change $x$ without changing the length of the interval), we assume in the following that
\be\label{norm}
\rho_0=1, \quad u_0,e_0, u_1, e_1 >0.
\ee

\subsection{Profile equations and formulation as mapping problem}\label{s:profile}

Our main interest is the study of steady solutions, i.e. solutions of
\ba\label{ns}
(\rho u)_x &=0\,,\\
(\rho u^2+\Gamma \rho e)_x &=\alpha u_{xx}, \\
\left( \rho u \left( e+\frac{u^2}{2} \right)+\Gamma \rho eu \right)_x &=\nu  e_{xx}+\left(\alpha  uu_x \right)_x\,
\ea
together with \eqref{BC}. In order to find these steady solutions, we use a shooting method. Integrating \eqref{ns} from $0$ to $x$ and rearranging using \eqref{norm}, there exists constants of integration $c=(c_{1},c_{2})$ to be determined so that we obtain similarly as in \cite{HLyZ2} the profile ODE
\ba\label{prof}
\frac{\alpha}{u_0} u' &= c_1 + u + \Gamma \frac{e}{u},\\
\frac{\nu}{u_0}e'&= c_2 -c_1 u - \frac{1}{2} u^2 + e,
\ea
together with $\rho=\frac{u_0}{u}$ and with the initial data
\ba\label{IC_prof}
u(0)=u_{0}>0,\\
e(0)=e_{0}>0.
\ea
In this setting
\be\label{link_c_derivative_at_0}
c_{1}=\frac{\alpha}{u_{0}} u'(0) - u_{0} - \Gamma \frac{e_{0}}{u_{0}} \text{ , } c_{2}=\frac{\nu}{u_{0}} e'(0) + \alpha u'(0) - e_{0} - \frac{1}{2} u_{0}^{2} - \Gamma e_{0},
\ee
where $u_{0}$ and $e_{0}$ are given and $(u'(0),e'(0))$ has to be determined in order to satisfies $(u(1),e(1))=(u_{1},e_{1})$.

\medskip

The domain of the ODE is the set
$$
 \left\{ (u,e) \in \R^{2} \text{ , } u > 0 \right\} 
$$
for which the right hand side of \eqref{prof} is well-defined and from which we can reconstruct $\rho$. 
Indeed, we remark that $u>0$ is imposed by $\rho u=\const$ and $\rho>0$. Hence we may ignore the variable $\rho$ in the following. Note also that the physical solutions are the ones for which $e$ is also positive. 

For a fixed choice of left data $(\rho_0, u_0, e_0)$ (meaning, by our previous normalization, 
just a fixed choice of $u_0$ and $e_0$), we define now the ``Cauchy-to-boundary value'' mapping
\be\label{Psi}
\Psi: (c_1,c_2) \to (u,e)(1),
\ee
where $(u,e)$ denotes the maximal solution of \eqref{prof}-\eqref{IC_prof}  for the given value of $c=(c_1,c_2)$.
Evidently, solutions of \eqref{BC}-\eqref{ns} thus correspond to solutions of
the mapping problem
\be\label{fund}
\Psi(c)=(u_1,e_1).
\ee

\section{The feasible set}\label{s:feasible}
In \eqref{Psi}, we did not specify the domain of $c$.
It is indeed our first order of business to determine it.
For a fixed choice of left data $(u_0, e_0)$,
we define the feasible set $\mathcal{C}_{u_{0},e_{0}}$ as the set of all $c$ for which \eqref{prof}-\eqref{IC_prof} has
a continuous solution $(u,e)$ on $[0,1]$ where $u$ and $e$ are both positive on $[0,1]$. Note that $\mathcal{C}_{u_{0},e_{0}}$ is not empty since $(-u_{0} - \Gamma \frac{e_{0}}{u_{0}}, - (1+\Gamma) e_{0} - \frac{1}{2} u_{0}^{2}) \in \mathcal{C}_{u_{0},e_{0}}$ (that corresponds to the constant solution of problem \eqref{prof}-\eqref{IC_prof}). Then, we have the following crucial observation.

\bpr\label{Cboundary}
The set $\mathcal{C}_{u_{0},e_{0}}$ is open and its boundary consists of $c$ for which there exists continuous functions $(u,e)$ on $[0,1]$, solution of
 Problem \eqref{prof}-\eqref{IC_prof} on $[0,1)$, with $u,e$ both positive on $[0,1)$ and such that $e(1)=0$.
\epr

Before proving Proposition \ref{Cboundary}, we establish a preliminary result.

\bl\label{small}
Let $c=(c_{1},c_{2}) \in \R^{2}$. Let $x_{*} \in (0,1]$ such that a solution $(u,e)$ of \eqref{prof}-\eqref{IC_prof} is defined on $[0,x_{*})$. We have the following statements :

(i) If $e > 0$ on $[0,x_{*})$, then $(u,e)$ is bounded on $[0,x_{*})$ uniformly with respect $x_{*}$ and one can extend continuously $(u,e)$ to $x_{*}$.

(ii) If there exists a constant $\tilde e>0$, $e \geq \tilde{e}>0$ on $[0,x_{*})$, then there exists a constant $\tilde{u}>0$, $u \geq \tilde{u}$ on $[0,x_{*})$.

%

(iii) If $(u,e)\to 0$ simultaneously as $x\to x_*$ with $(u,e)$ both positive on $[0,x_{*})$, 
then $c_1<0$ and $c_2<0$.

(iv) Assume $c_1<0$, $c_2<0$, $(u,e)$ is a solution of \eqref{prof}-\eqref{IC_prof} on $[0,x_{*}]$ and $u,e>0$ on $[0,x_*]$. For any $\eps>0$,
there exists $\delta>0$ depending only on $c_{1},c_{2},\eps$ such that if $u(x_*), e(x_*)\leq \delta$, there exists $\tilde x\in [x_*,x_*+ \eps]$ such that
$(u,e)$ extends continuously as a solution of \eqref{prof}-\eqref{IC_prof} on $[0, \tilde x)$ with $u,e>0$ on $[0,\tilde x)$ and $e(x)\to 0$ as $x\to \tilde x$.
\el

\br\label{remark_small}
As we will see in the proof of (iv), one can prove that if $c_1<0$, $c_2<0$ and $(u,e)$ is a solution of \eqref{prof}-\eqref{IC_prof} on $[0,1)$ with $u$ and $e$ both positive on $[0,1)$, there exists a constant $M>0$ depending only on $c_{1},c_{2},\Gamma$, such that for any $\delta>0$ small enough, for any $x_{\ast} \in (0,1)$, if $0< u(x_*), e(x_*)\leq \delta$, then $0< e \leq \delta$ and $0< u \leq M\delta$ on $[x_{*},1)$.
\er

\begin{proof}[Proof of Lemma \ref{small}]
(i) Since $e > 0$ on $[0,x_{*})$, we get
\begin{align*}
\frac{1}{2} \left( \frac{\alpha}{u_0} u^2+ \frac{\nu}{u_0}e^2 \right)' &= c_1 u + u^2 + \Gamma e + c_2 e - c_1 eu - eu^2/2 + e^2\\
&\leq c_1 u + u^2 + \Gamma e + c_2 e - c_1 eu + e^2\\
&\leq \frac12 \left( c_1^2 + u^2 + 2u^2 + \Gamma^2 + e^2 + c_2^2 + e^2 + c_1^2 e^2 + u^2 + e^2 \right)\\
&\leq A  \left( \frac{\alpha}{u_0} u^2+ \frac{\nu}{u_0}e^2 \right) + B
\end{align*}
for some constants $A,B>0$ depending only on $c_{1}$, $c_{2}$, $\Gamma$, $\alpha$, $\nu$, $u_{0}$. Hence, $|(u,e)|$ grows at most exponentially, in particular remaining bounded on $[0,x_{*})$. Furthermore, $u$ and $e$ can be continuously extended to $x_{\ast}$ since $(u^2)'$ and $(e^2)'$ are bounded and then integrable on $[0,x_{*})$. 
\\
(ii) The term $\Gamma \frac{e}{u}$ in the $u$-equation serves as a barrier meaning that there exists $u_{\ast}>0$ such that for any $u \in (0,u_{\ast}]$, $c_1 + u + \Gamma \frac{\tilde e}{u} \geq 0$. 
\\
(iii) Evidently, $c_1<0$, or else $u'>0$ for $u$, $e>0$, contradicting the assumed convergence to $0$.
Then, for $u>0$ sufficiently small, this implies that $-c_1 u-\frac{1}{2} u^2>0$ and hence
$\frac{\nu}{u_0} e'>c_2 +e$. Therefore, $c_2<0$ or else $e'>0$ for $e>0$ and $u>0$ sufficiently small, again contradicting convergence.
This proves (iii).

(iv) We assume now that $c_1$, $c_2<0$, that $(u,e)$ are both positive on $[0,x_*]$ and that $u(x_*), e(x_*)\leq \delta$. By point (i), one can extend $(u,e)$ as a solution of \eqref{prof}-\eqref{IC_prof} on a interval that strictly contains $[0,x_{\ast}]$. We introduce
\[
\tilde x = \sup \{ x \leq x_{\ast} + \eps \text{ , } (u,e) \text{ extends as a solution of \eqref{prof}-\eqref{IC_prof} and are positive on } [0,x) \}
\]
and we keep the notation $(u,e)$ for the solution on $[0,\tilde x)$. We first note that so long as $u$ and $e$ remain less than  $\frac{c_{2}}{c_{1}-1}$ on $[x_{\ast} , \tilde x)$
we have $\frac{\nu}{u_{0}} e' < - \frac{1}{2} u^2 \leq 0$ and thus $e$ is decreasing on $[x_{\ast} , \tilde x)$.
Next, based on the $u$-equation, several situations can happen around $x \in [x_{\ast} , \tilde x)$ :

\noindent (a) If $u(x)>\frac{-c_{1}-\sqrt{c_{1}^{2} - 4 \Gamma e(x)}}{2}$, $u'(x)<0$ and $u$ is decreasing around $x$.

\noindent (b) If $u(x)=\frac{-c_{1}-\sqrt{c_{1}^{2} - 4 \Gamma e(x)}}{2}$, $u'(x)=0$, $u''(x)= \Gamma \frac{e'(x)}{u(x)} <0$, $u$ is decreasing around $x$.

\noindent (c) If $u(x)<\frac{-c_{1}-\sqrt{c_{1}^{2} - 4 \Gamma e(x)}}{2}$, then $u<\frac{-c_{1}-\sqrt{c_{1}^{2} - 4 \Gamma e}}{2} \leq 2 \Gamma \frac{|e|}{|c_{1}|}$ around $x$.

\noindent Therefore, for $\delta$ small enough, $u \leq \max \left(1, \frac{2 \Gamma}{|c_{1}|} \right) \delta$ and $\frac{\nu}{u_0} e'<\frac{c_2}{2}$ on $[x_{\ast} , \tilde x)$ and $e$ goes to zero as $x \to \tilde x$ with
$|\tilde x-x_*|\leq \frac{2 \nu}{u_{0} |c_{2}|} \delta$. This proves assertion (iv).
\end{proof}

Thanks to this lemma we can assert that
$$
\mathcal{C}_{u_{0},e_{0}} = \left\{ c \in \R^{2} \text{ , where } e >0 \text{ on } [0,1] \text{ , } (u,e) \text{ the maximal solution of \eqref{prof}-\eqref{IC_prof}}   \right\}.
$$

\begin{proof}[Proof of Proposition \ref{Cboundary}]
For $c=(c_{1},c_{2})$, we denote by $(u,e)$ the maximal solution of Problem \eqref{prof}-\eqref{IC_prof}. Note that the following map is locally Lipschitz
\be\label{mapODE}
\Phi:(u,e) \in \left\{ (u,e) \in \R^{2} \text{ , } u > 0 \right\} \mapsto \left( c_{1} + u + \frac{\Gamma e}{u} \text{  ,  } c_{2} - c_{1}u - \frac{1}{2} u^2 + e \right).
\ee
If $c \in \mathcal{C}_{u_{0},e_{0}}$, then $u$ and $e$ are defined and positive on $[0,1]$ and by continuous dependence on parameters of solutions of an ODE, $c$ lies in the interior of $\mathcal{C}_{u_{0},e_{0}}$. In particular $\mathcal{C}_{u_{0},e_{0}}$ is open.
\\
We assume in the following that $c \in \mathcal{C}_{u_{0},e_{0}}^c$. Thanks to Lemma \ref{small}(i)-(ii), there exists $x_{*} \in (0,1]$ such that $u$ and $e$ are defined and continuous on $[0,x_*)$, $u,e>0$ on $[0,x_{*})$ and $(u,e)$ can be extended to $x_{*}$ with $e(x_*)=0$ and $u(x_{*})\geq 0$. Our goal is to show that $c \in \partial \mathcal{C}_{u_{0},e_{0}}$ if and only if $x_{\ast}=1$. Three different situations can then occur.
\\
Case (i) : $x_{\ast}=1$ and then $e(1)=0$.
\\
Case (ii) : $u(x_{*})>0$ and $x_{*} <1$. Then $(u,e)$ is defined on an interval that strictly contains $[0,x_{\ast}]$ and $e'(x_{\ast}) \leq 0$. In that case let us show that $e$ must actually cross $0$ and must become negative as $x$ crosses $x_*$. Since $u(x_{*})>0$, $(u,e)$ is defined on an interval that strictly contains $[0,x_{*}]$ and $e'(x_{*}) \leq 0$. Several subcases occur.

Subcase (ii)(a) : If $e'(x_{*}) < 0$, $e$ crosses $0$ and becomes negative as $x$ crosses $x_*$.

Subcase (ii)(b) : If $e'(x_*) = 0$ and $u(x_*) \neq -c_1$, then
$$
\frac{\nu}{u_0}e''(x_*)= -(c_1+u(x_*)) u'(x_*)=- \frac{u_{0}}{\alpha} (c_1+u(x_*))^2 < 0
$$
and $e$ crosses $0$ and becomes negative as $x$ crosses $x_*$.

Subcase (ii)(c) : If $e'(x_*) = 0$ and $u(x_*) = -c_1$, repeated differentiation shows that derivatives of $e$ and $u$ at $x_{*}$ vanish to all orders. By analyticity of solutions of an analytic ODE (note that $u>0$), $e \equiv 0$ and $u \equiv -c_{1}$, contradicting $e(0)=e_{0}>0$ so that this subcase can not occur.
\\
With such a fact in hand, there exists $\eps>0$ small enough such that $(u,e)$ is defined on $[0,x_{*}+\eps]$, $e$ negative on $]x_{*},x_{*}+\eps]$ and $u$ positive on $[0,x_{*}+\eps]$. By continuous dependence on parameters that $c$ lies in the interior  of $\mathcal{C}_{u_{0},e_{0}}^c$ and $c \not \in \partial \mathcal{C}_{u_{0},e_{0}}$.
\\ 
Case (iii) : $u(x_{*})=0$  and $x_{*} <1$.
In this case,  Lemma \ref{small}(iii) shows that $c_{1}$, $c_{2}<0$.  Let $(\tilde{c}_{1},\tilde{c}_{2})$ close enough to $(c_{1},c_{2})$ and denote by $(\tilde{u},\tilde{e})$ the maximal solution of Problem \eqref{prof}-\eqref{IC_prof} associated to $(\tilde{c}_{1},\tilde{c}_{2})$. We then take a $\delta$ associated to $\eps =1-x_{\ast}$ in Lemma \ref{small}(iv) that works for any $(\tilde{c}_{1},\tilde{c}_{2})$ close enough to $(c_{1},c_{2})$. By continuity of $u$ and $e$, there exists a number $\mu>0$ small enough such that $0<u(x_{*}-\mu),e(x_{*}-\mu) \leq \frac{\delta}{2}$. Then, by continuous dependence on parameters, for any $(\tilde{c}_{1},\tilde{c}_{2})$ close enough to $(c_{1},c_{2})$, $(\tilde{u},\tilde{e})$ is defined on $[0,x_{*}-\mu]$ and $0<\tilde{u}(x_{*}-\mu),\tilde{e}(x_{*}-\mu) \leq \delta$. Lemma \ref{small}(iv) shows that there exists $\tilde{x} \in [x_{\ast}-\mu,1-\mu]$, $\tilde{e}(\tilde{x})=0$. In particular, in that case $c \not \in \partial \mathcal{C}_{u_{0},e_{0}}$.
\end{proof}

We can now show that $\Psi$ defined in \eqref{Psi} is continuous.

\bpr\label{Psi_C0}
The map $\Psi$ is continuous on $\mathcal{C}_{u_{0},e_{0}}$ and can be extended to $\overline{\mathcal{C}_{u_{0},e_{0}}}$ as a continuous map denoted again $\Psi$.
\epr

\begin{proof}
The fact that $\Psi$ is continuous on $\mathcal{C}_{u_{0},e_{0}}$  follows from continuous dependence on parameters of solutions of an ODE (and the fact that the map $\Phi$ defined in \eqref{mapODE} is locally Lipschitz). We consider now $c=(c_{1},c_{2}) \in \partial \mathcal{C}_{u_{0},e_{0}}$. Proposition \ref{Cboundary} shows the maximal solution $(u,e)$ of \eqref{prof}-\eqref{IC_prof} is defined and continuous on $[0,1)$ and can be extended continuously to $1$ with $e(1)=0$ and $u(1) \geq 0$. Therefore, we can define $\Psi(c)=(u(1),0)$. If $u(1)>0$, $(u,e)$ is defined on an interval that strictly contains $[0,1]$ and by continuous dependence on parameters, $\Psi$ is continuous at $c$. We now have to deal with the case  $u(1)=0$. Lemma \ref{small}(iii) shows that $c_{1}$, $c_{2}<0$. Consider $\eps>0$. Let $\tilde{c} \in \overline{\mathcal{C}_{u_{0},e_{0}}}$ close enough to $c$ and denote by $(\tilde{u},\tilde{e})$ the maximal solution of Problem \eqref{prof}-\eqref{IC_prof} associated to $\tilde{c}$. By continuity of $(u,e)$ there exists $x_{*} \in (0,1)$ such that $0<u(x_{*}),e(x_{*}) \leq \frac{\eps}{2}$. Then, by continuous dependence on parameters, for any $\tilde{c}$ close enough to $c$, we have $0<\tilde{u}(x_{*}),\tilde{e}(x_{*}) \leq \eps$. Finally, by taking $\eps$ small enough, Remark \ref{remark_small} and the fact that $\tilde{u}$ is continuous at $1$ give $0<\tilde{u}(1) \leq M \eps$ (where $M$ depends only on $c$ and $\Gamma$). Hence, $\Psi$ is continuous at $c$.
\end{proof}

\section{Existence}\label{s:existence}
We are now ready to study existence. We first show that $\Psi$ is ``proper'' in the following sense.

\bpr\label{properprop}
Assume that $u_0>0,e_0>0$ are fixed. Let $c=(c_{1},c_{2}) \in \mathbb{R}^{2}$, such that $|c| \gg 1$ and denote by $(u,e)$ the maximal solution of \eqref{prof}-\eqref{IC_prof}. Then, if $c \in \mathcal{C}_{u_{0},e_{0}}$, either $u(1)\gg 1$, $e(1)\gg 1$ or $0< u(1)\ll 1$.
\epr

\begin{proof}
Several situations can happen.

Case (i) ($c_1 \gg 1$). If $c \in \mathcal{C}_{u_{0},e_{0}}$, $\frac{\alpha}{u_0} u' \geq c_1 + u_{0}$ and $u(1) \geq \frac{u_0}{\alpha} (c_{1} + u_{0}) + u_{0} \gg 1$. 

Case (ii) ($c_2\ll -1$). We consider the energy $y=\frac{\alpha}{2 u_{0}} u^2 + \frac{\nu}{u_{0}} e$. Then,
$$
y' = c_{2} + (\Gamma + 1) e + \frac{u^2}{2} \leq c_{2} + My
$$
where $M$ is a positive constant depending only on $\Gamma, u_{0},\alpha, \nu$. Then, for any $x$ in the domain of definition of $y$,
\[
y(x) \leq y(0) e^{Mx} + \frac{c_{2}}{M}(e^{Mx}-1),
\]
so that if $c_{2} < -\frac{M}{e^{M}-1} y(0) e^{M}$, $e$ must vanish at a point between $0$ and $1$ and $c \not \in \mathcal{C}_{u_{0},e_{0}}$.

Case (iii) ($c_2\gg 1$).
Let $c \in \mathcal{C}_{u_{0},e_{0}}$. We consider again the energy $y=\frac{\alpha}{2 u_{0}} u^2 + \frac{\nu}{u_{0}} e$. Then,
$$
y' = c_{2} + (\Gamma + 1) e + \frac{u^2}{2} \geq c_{2} + my
$$
where $m$ is a constant depending only on $\Gamma,u_{0},\alpha,\nu$. Therefore either $e(1) \gg 1$ or $u(1) \gg 1$.

Case (iv) ($c_1 \ll -1$ and $c_{2} \leq \sqrt{-c_{1}}$)
Let $c \in \mathcal{C}_{u_{0},e_{0}}$. Using again $y=\frac{\alpha}{2 u_{0}} u^2 + \frac{\nu}{u_{0}} e$ and following case (ii), we get $y' \leq -c_{1} + My$ on $[0,1]$. Therefore, for $c_{1} \ll -1$, there exists a constant $B$ depending only on $u_{0},e_{0},\alpha,\nu,\Gamma$ such that $e \leq B \sqrt{-c_{1}}$ on $[0,1]$. Using this fact on the u-equation of System \eqref{prof}, we get  
$$
\frac{\alpha}{u_{0}} u' \leq c_{1} + u + \frac{\Gamma B \sqrt{-c_{1}}}{u}.
$$
Note that for $x \in [0,1]$, if $u(x) \in \left[ - \frac{c_{1}}{4} - \frac{1}{2} \sqrt{\frac{1}{4} c_{1}^{2} - 4 \Gamma B \sqrt{-c_{1}}} , - \frac{c_{1}}{4} + \frac{1}{2} \sqrt{\frac{1}{4} c_{1}^{2} - 4 \Gamma B \sqrt{-c_{1}}} \right]$, then $\frac{\alpha}{u_{0}} u'(x) \leq \frac{1}{2} c_{1}$ and denote $b=- \frac{c_{1}}{4} - \frac{1}{2} \sqrt{\frac{1}{4} c_{1}^{2} - 4 \Gamma B \sqrt{-c_{1}}}$. We then notice that, for $c \in \mathcal{C}_{u_{0},e_{0}}$ and $c_{1} \ll -1$, $u$ rapidly goes under $b$ and stays under $b$. Therefore $u(1) \leq b  \leq \frac{4 \Gamma B}{\sqrt{-c_{1}}}$ and $u(1) \ll 1$.
\end{proof}

\begin{remark}
We proved in the previous proposition that there exists a constant $A>0$ depending only on $u_{0},e_{0},\alpha,\nu,\Gamma$ such that for any $c_{2} \leq -A$, $c \not \in \mathcal{C}_{u_{0},e_{0}}$ (see case(ii)). Note also that if $\frac{\nu}{u_{0}} e_{0}+c_{2}-c_{1} u_{0} - \frac{1}{2} u_{0}^{2} + e_{0} \leq 0$ and $c_{1} + u_{0} \geq 0$, then $c \notin \mathcal{C}_{u_{0},e_{0}}$. Indeed, in this case, $u$ increases and $-c_{1} u - \frac{1}{2}u^2 \leq -c_{1} u_{0} - \frac{1}{2}u_{0}^2$ so that $e$ is decreasing, $e' \leq - e_{0}$ and then $e$ crosses $0$ in $(0,1]$.
\end{remark}

The previous proposition is not empty in the sense that $\mathcal{C}_{u_{0},e_{0}}$ is not bounded. 

\bl\label{C_unbounded}
Assume $u_0>0,e_0>0$ are fixed. There exists a positive number $A$ depending only on $\Gamma,u_{0},e_{0},\alpha,\nu$ such that if $c_{1} + u_{0} + \frac{\Gamma A (|c_{2}|+1)}{u_{0}} < 0$ and $c_{2} + e_{0} > 0$, then $(c_{1},c_{2}) \in \mathcal{C}_{u_{0},e_{0}}$.
\el

\begin{proof}
For $c=(c_{1},c_{2})$, we denote by $(u,e)$ the maximal solution of Problem \eqref{prof}-\eqref{IC_prof} and by $I$ its interval of definition. Following case (ii) in the previous proposition there exists a constant $A>0$ depending only on $\Gamma,u_{0},e_{0},\alpha,\nu$ such that $e \leq A (|c_{2}|+1)$ on $[0,1] \cap I$. Then we note that since $c_{1} < c_{1} + u_{0} <0$, the map $y \in [0,u_{0}] \mapsto -c_{1} y - \tfrac12 y^{2}$ is nonnegative so that for any $x \in [0,1] \cap I$ such that $0<u(x) \leq u_{0}$, we have $\frac{\nu}{u_0}e'(x) \geq c_{2} + e(x)$. Note also that for $x \in [0,1] \cap I$ such that $u(x) = u_{0}$,
$$
\frac{\alpha}{u_0} u'(x) \leq c_1 + u_{0} + \Gamma \frac{ A (|c_{2}|+1)}{u_{0}}  < 0.
$$
In particular $u'(0)<0$, $e'(0) > c_{2} + e_{0} > 0$ and $\{ x \in (0,1] \cap I \text{ , } u(x)=u_{0} \}$ is empty. By Lemma \ref{small}(i)-(ii), $[0,1] \subset I$, $e$ is increasing and $0 < u \leq u_{0}$ on $[0,1]$so that $c \in \mathcal{C}_{u_{0},e_{0}}$.
\end{proof}

We now define for $\eps>0$, $E_{\eps} = \{(x,y) \in \mathbb{R}^2 \text{ , } \eps < x,y < \frac{1}{\eps} \}$ and $\Omega_{\eps}=\Psi^{-1}(E_{\eps})$. By continuity of $\Psi$ (Proposition \ref{Psi_C0}) and Proposition \ref{properprop}, $\Omega_{\eps}$ is open, bounded and $\overline{\Omega_{\eps}} \subset \mathcal{C}_{u_{0},e_{0}}$.  
We denote by $\Psi_{\eps}$ the restriction of $\Psi$ to $\Omega_{\eps}$,
and by $d(\Psi_{\eps},\Omega_{\eps},(e_{1},u_{1}))$ the Brouwer degree \cite{Brouwer1911,Mil1965,Hi1976,Pra2006,DincaMawhin2021}
of $\Psi_{\eps}$ in $\Omega_{\eps}$ with respect to the target $(e_{1},u_{1}))$. 

Recall, for {\it regular values} $(e_1,u_1)$, defined as values for which $\Psi_\eps$ is differentiable and full rank on
$\Psi_\eps^{-1}(e_1,u_1)$, 
\be\label{regdegree}
d(\Psi_{\eps},\Omega_{\eps},(e_{1},u_{1})):=\sum_{c\in \Psi_\eps^{-1}(e_1,u_1)} \sgn \det d\Psi_\eps(c), 
\ee
that is, the degree counts roots with sign depending on orientation.
For arbitrary (not necessarily regular) values $(e_1,u_1)$ for which $\Psi_\eps^{-1}(u_1,e_1)=\emptyset$,
$d(\Psi_{\eps},\Omega_{\eps},(e_{1},u_{1}))=0$.  Thus, {\it nonzero Brouwer degree implies existence of a solution}.
Finally, recall that $d(\Psi_{\eps},\Omega_{\eps},(e_{1},u_{1}))$ is homotopy invariant, so long as
$\Psi_\eps^{-1}(e_1,u_1)$ remains disjoint from 
$\partial \Omega_\eps$.
Typically, degree is evaluated at a regular value, then deduced for other values by homotopy invariance.

\bc\label{propcor}
Assume that $u_0>0$ and $ e_0>0$ are fixed. Let $u_{1}>0$, $e_{1}>0$. Then for $\eps>0$ small enough, $(u_{1},e_{1}) \not \in \Psi(\partial \Omega_{\eps})$ and the Brouwer degree $d(\Psi_{\eps},\Omega_{\eps},(e_{1},u_{1}))$ is independent of $(u_{1},e_{1})$ and $\eps$.
\ec

\begin{proof}
Let $u_{1}>0,e_{1}>0$. First, Proposition \ref{properprop} shows that $\Psi^{-1}(u_{1},e_{1})$ is bounded and included in the open set $\Omega_{\eps}$ for $\eps$ small enough. In particular, $d(\Psi_{\eps},\Omega_{\eps},(e_{1},u_{1}))$ is independent of $\eps$ small enough. Furthermore, we also get from Proposition \ref{properprop} that for any $t\in[0,1]$, $(1+(1-t)u_{1},1+(1-t)e_{1}) \not \in \Psi(\partial \Omega_{\eps})$ if $\eps$ is small enough. Hence, by homotopy invariance, $d(\Psi_{\eps},\Omega_{\eps},(e_{1},u_{1}))$ and $d(\Psi_{\eps},\Omega_{\eps},(1,1))$ are equal.
\end{proof}

At this point we make use of the fundamental property that gas dynamics has an associated convex entropy
$\eta(\rho, \rho u, \rho E)$ in the sense of \cite{La,KSh}, namely $\eta=-\rho S$, where $S(\rho^{-1} , e)$ 
is thermodynamic entropy; see \cite[\S 4]{KSh}.
That is, writing \eqref{nsint} as $ U_t + f(U)_x=(B(U)U_x)_x $,
where $U=(\rho, \rho u, \rho E)^t$, there hold: for any $U \in \R^{3}$, the Hessian matrix
$d^2\eta_{U}>0$; $d\eta_{U} \circ df_{U}=dq_{U}$ for some flux $q$; and $d^2\eta_{U} B(U)$ is symmetric positive semidefinite, hence
$\langle V,  d^2\eta_{U} B(U) V\rangle \leq 0$, with equality if and only if $d^2\eta_{U} B(U)V=0$, or equivalently $B(U)V=0$. 
Composing the equations on the left by $d\eta$, we have 
\be\label{entprelim}
\begin{aligned}
	\eta(U)_t + q(U)_x&= d\eta_{U} (B(U)U_x)_x= (d\eta_{U} (B(U) U_x))_x - \langle d^2\eta_{U} U_x, B(U) U_x\rangle\\
	&= (d\eta_{U} (B(U) U_x))_x - \langle U_x, d^2\eta_{U} B(U) U_x\rangle,
\end{aligned}	
	\ee
giving $\eta(U)_t + q(U)_x - (d\eta_{U} (B(U)U_x))_x\leq 0$ with equality if and only if $B(U)U_x=0$.
Using the definition $B(U)U_x= (\alpha u_{x}, \kappa T_{x}+\alpha uu_x)^t$ given by \eqref{nsint},
we find that $B(U)U_x=0$ is equivalent to $(u,T)_x=0$, and therefore to $(u,e)_x=0$.
Integrating the time-independent profile equation from $x=0$ to $1$, we thus obtain
\be\label{entineq}
\hbox{\rm $\big[ q(\hat U) - d\eta_{\hat U} (B(\hat U) \hat U_x) \big]_0^1 \leq 0$,
with equality if and only if $(\hat u,\hat e)\equiv \const$.}
\ee

\bl\label{entlem}
For steady gas dynamics on an interval, \eqref{nsint}--\eqref{BC},
with constant boundary conditions $(u_0,e_0)=(u_1,e_1)$, the unique global
solution is given by the constant solution $(\hat \rho, \hat u, \hat e)\equiv (\rho_0, u_0, e_0)$.
Equivalently $c_*= (- u_{0} - \Gamma e_{0}/u_{0},   -(1+\Gamma) e_{0} -  u_{0}^{2}/2 )$
is the unique global solution of $\Psi(c)=(u_1,e_1)$; moreover, it is nondegenerate, with $\sgn \det d\Psi(c_*)= +1$.
\el

\begin{proof}
	From $\hat \rho \hat u=: m \equiv \const$, we obtain $(\rho, u, e)(0)=(\rho, u, e)(1)$,
	or $U(0)=U(1)$ in the vectorial notation above.
	Note that $\eta$ may be modified by the addition of any affine function
	while preserving its properties as a convex entropy (by changing $q$ accordingly).
	Thus, by an appropriate affine shift, we may arrange that
	$\eta(U(0))=0$ and $d\eta_{U(0)}=0$, so that the 
	left-hand side vanishes in \eqref{entineq} (since $\eta(U(0))=\eta(U(1))$ and $d\eta_{U(0)}=d\eta_{U(1)}$), and therefore $(\hat u, \hat e)\equiv \const$.
	But, then, $\hat \rho=m/\hat u\equiv \const$ as well, and so $(\hat \rho, \hat u, \hat e)
	\equiv (\rho_0, u_0, e_0)$ as claimed.
	The computation of $d\Psi(c^*)$ amounts to integration of a $2\times 2$ 
	constant-coefficient linearized equations about this constant solution,
	hence may be carried out explicitly to find that $\sgn \det d\Psi(c_*)=+1$.
	We omit this calculation as we will show it in a simpler and more general way later on. See Subsection \ref{s:stabind}.
\end{proof}

\br\label{fwdrefrmk}
Lemma \ref{entlem} may be recognized as a particularly concrete instance of results stated for general
systems in \cite[Thm 2.10]{paper2} and \cite[Prop. 2.9]{paper2}.
In particular, $\sgn \det d\Psi(c_*)=+1$ is seen by abstract considerations
to hold for constant solutions of general symmetrizable systems, without explicit calculation.
We note for gas dynamics that the key identity \eqref{entprelim} may be obtained readily
from the thermodynamic relation $de=T\, dS + p\, dv$ defining $S$, where $V=1/\rho$, or $S_t= T^{-1}(e_t - pv_t)$,
together with \eqref{nsint}, without verifying convexity or symmetrizability, with no need to invoke general theory.
\er

\bt[Large-data existence]\label{maincor} 
For steady gas dynamics on an interval, \eqref{nsint}--\eqref{BC},
there is at least one steady solution for every choice of left and right data.
\et

\begin{proof}
Applying Corollary \ref{propcor}, we find that the Brouwer degree 
is independent of the target $(u_1, e_1)$. 
	Thus we may compute the degree at the constant data $(u_1, e_1) =(u_0, e_0)$.
	By Lemma \ref{entlem}, $\Psi^{-1}(u_1,e_1)$ consists of the single point
	$c_*= (- u_{0} - \Gamma e_{0}/u_{0},   -(1+\Gamma) e_{0} -  u_{0}^{2}/2 )$, at which
	$\sgn \det d\Psi(c_*)=+1$.  Thus, by \eqref{regdegree}, the degree at $(u_1,e_1)$ is $+1$.
	This implies that the Brouwer degree is $+ 1$ for all values of the target,
	implying existence of a solution.
\end{proof}

\section{Uniqueness}\label{s:section}

We next characterize uniqueness, by a global version of the local \cite[Lemma 3.10]{BFZ}.

\bpr\label{uniqueprop}
If $\gamma:=\det d\Psi(c)$ does not vanish on the feasible set $\mathcal{C}_{u_{0},e_{0}}$,  then solutions
of \eqref{fund} are globally unique for each choice of data $(\rho_0, u_0, e_0, u_1,e_1)$. 
If on the other hand $\gamma$ changes sign on the feasible set $\mathcal{C}_{u_{0},e_{0}}$, then even local uniqueness is violated;
in particular, there is at least one choice of data possessing multiple solutions.
\epr

\begin{proof}
Nonvanishing of $\gamma$ implies nonvanishing of the Jacobian determinant $\det d\Psi(c)$, 
	which implies $d\Psi(c)$ full rank and $\sgn\det d\Psi(c)\equiv +1$ for all $c$. In particular all 
	values are regular and it follows that the degree of $\Psi$ with respect to a target $(u_1, e_1)$
is equal to $+ n$, where $n$ is the number of solutions for that data. 
Since we have already shown that degree is identically equal to $+ 1$, this is a contradiction unless 
roots are unique i.e., $n=1$. This proves the first assertion. For the second assertion, just notice that uniqueness 
implies that degree is equal to the sign of $\gamma$ at the unique solution and therefore a change of sign in
$\gamma$ implies a change in degree. 
Hence, by contradiction, uniqueness is impossible when $\gamma$ changes sign. 
\end{proof}

\noindent{\bf Conclusion}: {\it Uniqueness or nonuniqueness hinges on nonvanishing of $\det d\Psi(\cdot)$ on $\mathcal{C}_{u_{0},e_{0}}$.}

\section{Spectral stability and the Evans function}\label{s:evans}
We can reduce Problem \eqref{nsint} to
\begin{align*}
&\rho_{t} + \left( \rho u \right)_{x} = 0 \,, \\
&\rho u_{t} + \rho u u_{x} + (\Gamma\rho e)_{x} = \alpha u_{xx} \,, \\
&\rho e_{t} + \rho u e_{x} + \Gamma \rho e u_{x} = \nu e_{xx} + \alpha  u_{x}^{2} \,,
\end{align*}
from which we obtain the eigenvalue problem around a steady state $(\hat{\rho} ,\hat{u},\hat{e})$
\ba\label{eval}
&\lambda \rho + \left(\hat{\rho} u + \hat{u} \rho \right)_{x} = 0 \,, \\
&\lambda \hat{\rho} u + \left(\hat{\rho} \hat{u} u +\Gamma \hat{\rho} e + \Gamma \hat{e} \rho \right)_{x} + \hat{u}_{x} \left(\hat{\rho} u + \hat{u} \rho \right) = \alpha  u_{xx} \,, \\
&\lambda \hat{\rho} e + \left(\hat{\rho} \hat{u} e \right)_{x} + \hat{e}_{x} \left(\hat{\rho} u + \hat{u} \rho \right) + \Gamma \hat{\rho} \hat{e} u_{x} + \Gamma \hat{u}_{x} \left(\hat{\rho} e + \hat{e} \rho \right) = \nu e_{xx} + 2\alpha  \hat{u}_{x} u_{x} \,,
\ea
with boundary conditions 
\be\label{linbc}
(\rho, u,e)(0)=0, \; (u,e)(1)=0.
\ee
Note that for $\lambda=0$ the previous system can be written in the alternative form
\ba\label{eval_lambdazero}
& \hat{\rho} u + \hat{u} \rho = 0 \,, \\
&\left( \hat{\rho} \hat{u} u +\Gamma \hat{\rho} e + \Gamma \hat{e} \rho \right)_{x} = \alpha u_{xx} \,, \\
&\left((1 + \Gamma) \hat{\rho} \hat{u} e +\hat{\rho} \hat{u}^{2} u \right)_{x} = \nu e_{xx} + \alpha \left( \hat{u} u_{x} + \hat{u}_{x} u  \right)_{x} \,.
\ea

We are using here the standard approach \cite{AGJ,GZ} of rewriting 
\eqref{eval} as a first-order system and a Cauchy problem.
Note that, after eliminating $\rho$, \eqref{eval} may be rewritten as a first-order system in $(u, e,u',e')$,
following the standard approach of \cite{AGJ,GZ}, with homogeneous data prescribed on $(u,e)$ at both ends.
The {\it Evans function} may thus be defined via a ``shooting'' construction, 
similarly as in \cite{R,SZ} for the half-line case, as
\be\label{evans}
D(\lambda):= \det \bp u_1(1) & u_2(1) \\ e_1(1) & e_2(1)\ep,
\ee
where $(\rho_j, u_j, e_j)$ are solutions of \eqref{eval} with initial conditions 
$$
(\rho_1, u_1, e_1,u_1',e_1')(0)=(0,0,0,1,0 ), \quad  (\rho_2, u_2, e_2,u_2',e_2')(0)=(0,0,0,0,1 ); 
$$
that is, as the Wronskian at $x=1$ of a basis of solutions satisfying the boundary conditions at $x=0$.
This Wronskian vanishes precisely when there exists a solution vanishing in $(e,u)$ at both $x=0,1$,
i.e., an eigenfunction.
Evidently, $D(\cdot)$ is analytic in $\lambda$ on all of $\C$, and real-valued for $\lambda$ in $\R$, 
with zeros corresponding to eigenvalues of 
the linearized operator about the associated steady solution.\footnote{
Indeed, as standard in Evans function theory, zeros correspond in both location and multiplicity
to eigenvalues of the linearized operator about the wave; see, e.g., \cite{AGJ,GZ,ZH} in the whole-line case.}

\medskip

{\bf Conclusion}: {\it Spectral stability is equivalent to nonvanishing of $D$ on $\{\Re \lambda \geq 0\}$.}

\subsection{The stability index}\label{s:stabind}
Clearly $D$ is real-valued for real $\lambda$.
It is readily seen (see, e.g. \cite{MelZ})
that $D(\lambda)\neq 0$ for $\lambda$ real and sufficiently large, hence we may define as in \cite{GZ}
the {\it Stability index}
\begin{equation}\label{mudef}
\mu:= \sgn D(0)\left( \lim_{\lambda \to +\infty_{real}}\sgn D(\lambda) \right) 
\end{equation}
as a nonvanishing multiple $\pm \sgn D(0)$ of $\sgn D(0)$.
Evidently, $\mu$ determines the parity of the number of roots of the Evans function with positive real part, or, equivalently
(since complex roots occur in conjugate pairs), the number of positive real roots, with $+1$ corresponding to ``even''
and $-1$ to ``odd''.
As such, it is often useful in obtaining {\it instability} information.

Moreover, we have the following key observation 
relating the low-frequency stability and the stability index information to transversality of the steady profile
solution of the standing-wave ODE.

\bl\label{ZSlem}
The zero-frequency limit $D(0)$ is equal to $\frac{\alpha \nu}{u_{0}^{2}}$ multiplied by the
Jacobian determinant $\det d\Psi(c)$ associated with problem \eqref{fund} evaluated at any root $c$;
in particular,
\be\label{sgns}
\sgn D(0)=\sgn \det d\Psi(c).
\ee
\el

\begin{proof}
The proof amounts to the observation that the operations of linearization and integration of the standing-wave
ODE commute.
Taking the variation of the profile equation \eqref{prof} with respect to $c$ gives
\ba\label{linprof}
&\frac{\alpha}{u_0} \dot u' = \dot c_1 + \dot u + \Gamma \left( \frac{\dot e}{\hat{u}} - \frac{\hat{e}}{\hat{u}^2} \dot u \right)
,\\
&\frac{\nu}{u_0} \dot e'= \dot c_2 -(\dot c_1 \hat{u}+ c_1 \dot u) - \hat{u}\dot u + \dot e,\\
&(\dot u, \dot e)(0)=(0,0),
\ea
where $\dot{ }$ denotes variation. Furthermore, we deduce from relations \eqref{link_c_derivative_at_0} that 
$$(\dot c_1,\dot c_2)= \left( \frac{\alpha}{u_{0}} \dot u'(0), \frac{\nu}{u_{0}} \dot e'(0) + \alpha \dot u'(0) \right).
$$
It is readily verified for $\lambda=0$ that the eigenvalue equations \eqref{eval_lambdazero} can be integrated from $0$ to $x$
to yield the same system \eqref{linprof} (note that $\hat{\rho} \hat{u}=u_{0}$). Therefore, keeping the notations of \eqref{evans}, for $(\dot c_1,\dot c_2)=(1,0)$, $(\dot u,\dot e)=\frac{u_{0}}{\alpha} (u_{1},e_{1})-\frac{u_{0}^{2}}{\nu} (u_{2},e_{2})$, whereas for $(\dot c_1,\dot c_2)=(0,1)$, $(\dot u,\dot e)=\frac{u_{0}}{\nu} (u_{2},e_{2})$. The result follows.
\end{proof}

\begin{remark}
The previous lemma gives us another way to compute $D(0)$. Considering the problem
\ba\label{eqD0}
&\frac{\alpha}{u_{0}} u' = d_{1} + \left(1 - \Gamma \frac{\hat{e}}{\hat{u}^{2}} \right) u + \frac{\Gamma}{\hat{u}} e  \,, \\
&\frac{\nu}{u_{0}} e' = d_{2} - d_{1} \hat{u} -\frac{\alpha}{u_{0}} \hat{u}' u + e + \Gamma \frac{\hat{e}}{\hat{u}} u,\\
&u(0)=0 , e(0)=0,
\ea
we have
$$
D(0)= \det \bp u_1(1) & u_2(1) \\ e_1(1) & e_2(1)\ep
$$
where $(u_{1},e_{1})$ solves \eqref{eqD0} for $(d_{1},d_{2}) = \left( \frac{\alpha}{u_{0}}-u_{0}-\Gamma \frac{e_{0}}{u_{0}},\alpha-e_{0}-\frac{1}{2} u_{0}^{2}-\Gamma e_{0} \right)$ and $(u_{2},e_{2})$ solves \eqref{eqD0} for $(d_{1},d_{2}) = \left(-u_{0}-\Gamma \frac{e_{0}}{u_{0}},\frac{\nu}{u_{0}}-e_{0}-\frac{1}{2} u_{0}^{2}-\Gamma e_{0} \right)$ $ ($see \eqref{link_c_derivative_at_0} for the link between $(c_{1},c_{2})$ and $(u'(0),e'(0))$.
\end{remark}

\br\label{zsrmk}
Lemma \ref{ZSlem} is analogous to the Zumbrun-Serre/Rousset lemmas of \cite{ZS,R} in the whole- and half-line case,
which say $D(\lambda)\sim \gamma \delta(\lambda)$ for $|\lambda | \ll 1$, 
where $\gamma$ is a Wronskian encoding transversality of the associated standing-wave ODE and $\delta$ is a Lopatinski
determinant for the inviscid stability problem (here trivially nonvanishing).
\er

{\bf Conclusion}: Both Brouwer degree $\gamma=\sgn \det d\Psi ( \cdot)$ and stability index $\mu$ are determined by $\sgn (D(0))$, hence (by Proposition \ref{uniqueprop}
and the discussion just above) uniqueness and topological stability information may be obtained by evaluation of $D(0)$ on the feasible set $\mathcal{C}_{u_{0},e_{0}}$.
In particular, differently from the cases of the whole- or half-line (see, e.g., the discussion of \cite[\S 6.2]{Z2}),
{\it changes in stability/Morse index associated 
with passage of a single eigenvalue through $\lambda=0$ are necessarily associated with bifurcation/nonuniqueness.}

\section{Numerical investigations}\label{s:num}

For simple gases, the ratio $\frac{\nu}{\alpha}$ follows closely to the prediction
\be\label{stat}
\frac{\nu}{\alpha} 
= \frac{27 \Gamma +12}{16}
\ee
of statistical mechanics \cite{HLyZ1,HLyZ2}.\footnote{
In the notation of \cite{HLyZ2},  $\alpha=2 \mu + \eta= \frac{4}{3} \mu$, 
$\gamma=\Gamma+1$, and
$\frac{\kappa}{c_v \mu}= \frac{9\gamma -5}{4}$,
giving the result.
}
In our numerics, we will assume, further, \eqref{stat}.

\subsection{Feasible set}\label{s:numfeas} 
\begin{figure}[htbp]
 \begin{center}
\begin{equation*}
\begin{array}{lcr}
(a) \includegraphics[scale=0.21]{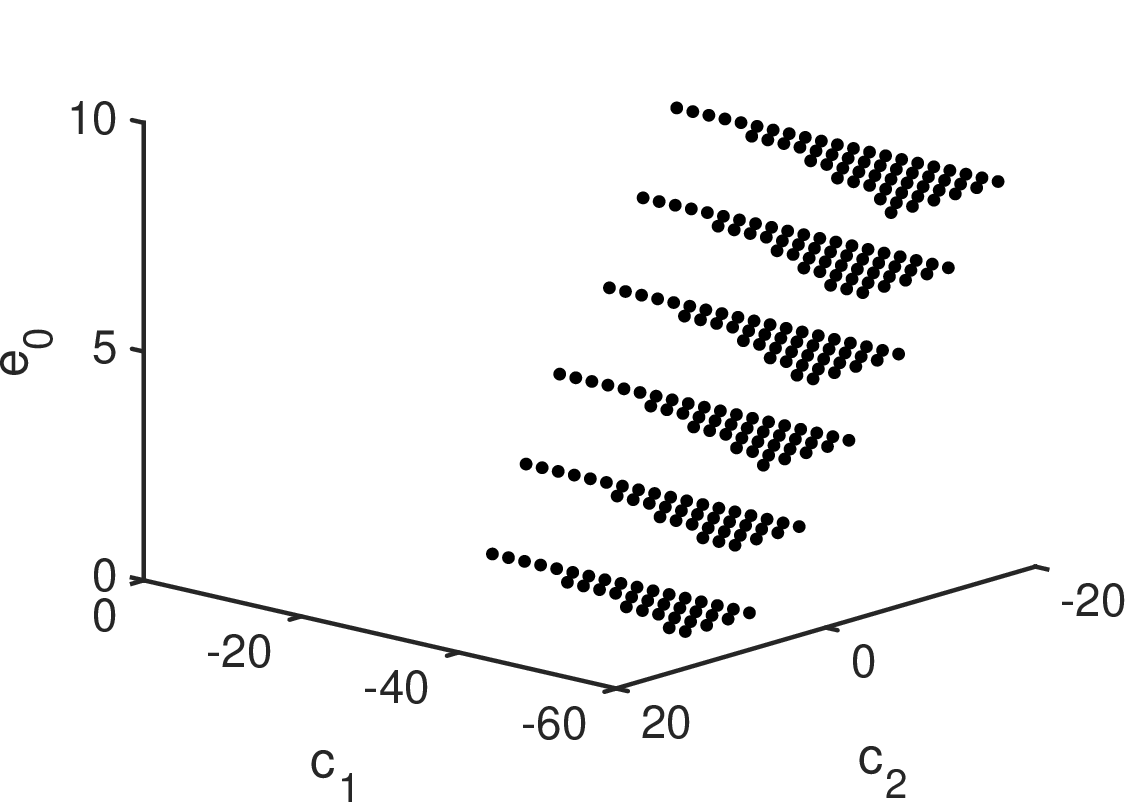}& (b) \includegraphics[scale=0.21]{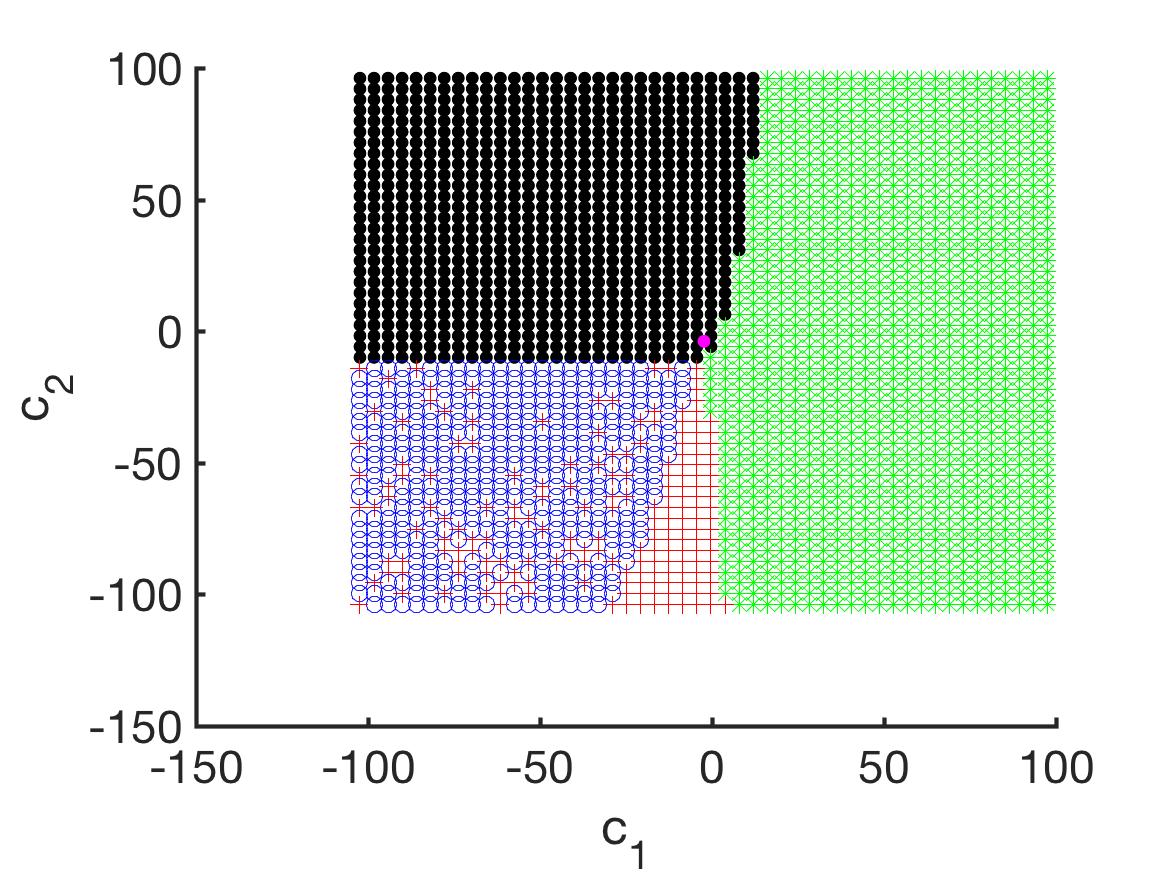}& (c) \includegraphics[scale=0.21]{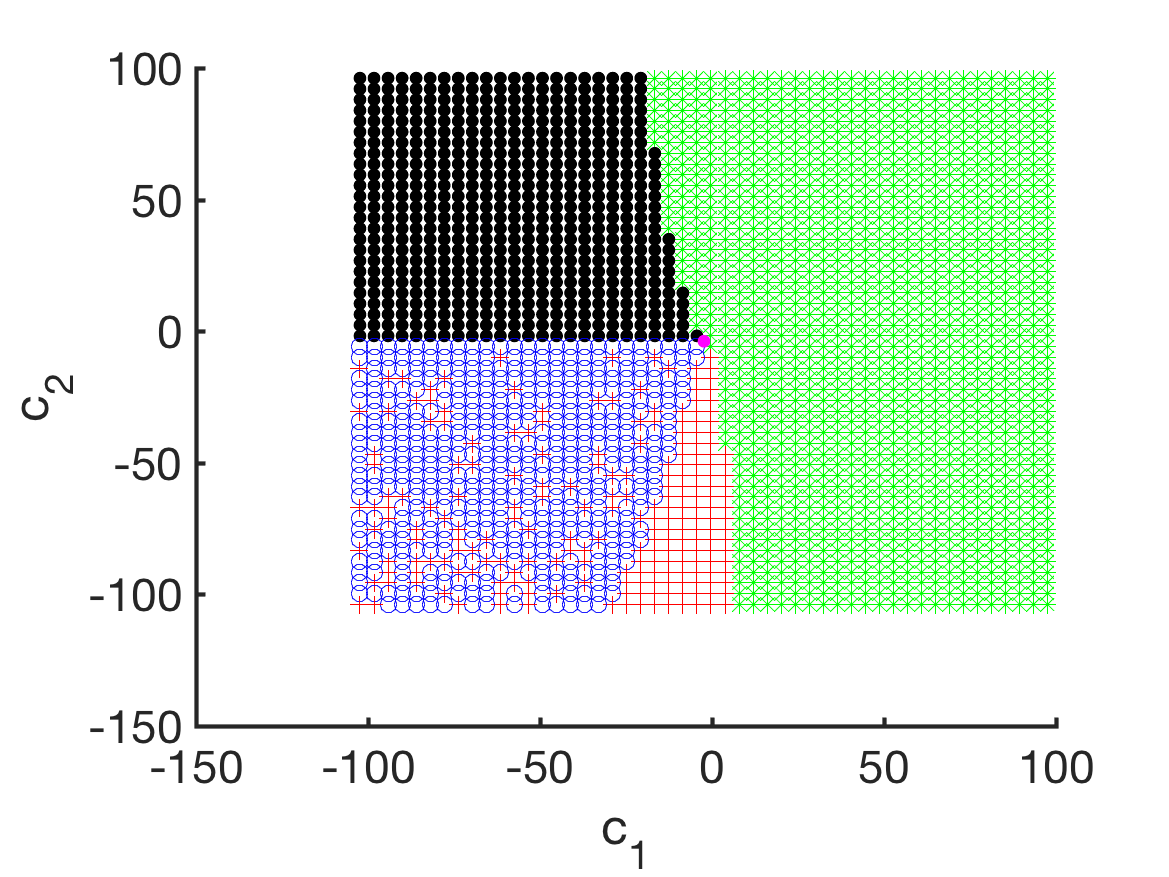}
\end{array}
\end{equation*}
\end{center}
\caption{(a) Plot of the feasible set as $e_0$ varies when $\Gamma = 1$, $\alpha = 0.1$,  and $\nu = 0.244$. (b) Plot of the feasible set with black dots, the set where $u$ goes negative on $[0,1]$ with blue circles, the set where $e$ goes negative on $[0,1]$ with green stars, and the set where there is finite time blowup on $[0,1]$ with red $+$ signs for  $\alpha = 2$, $\nu = 3.75$, $\Gamma = 2/3$, $e_0 = 2$. (c) Plot of the feasible set with black dots, the set where $u$ goes negative on $[0,1]$ with blue circles, the set where $e$ goes negative on $[0,1]$ with green stars, and the set where there is finite time blowup on $[0,1]$ with red $+$ signs for  $\alpha = 0.2$, $\nu = 1$, $\Gamma = 2/3$, $e_0 = 2$. A bold magenta dot marks the constant solution on plots (b) and (c). 
}
\label{fig:feasible}
\end{figure}

For our numerical studies, we rescale equation \eqref{nsint} through the following change of coordinates, $\rho = \rho_0 \bar \rho$, $u = u_0\bar u$, $e = u_0^2 \bar e$, $t = \frac{\bar t}{u_0}$,  $\bar \alpha := \frac{\alpha}{\rho_{0} u_0}$, $\bar \nu := \frac{\nu}{\rho_{0} u_0}$, which allows us to always fix $\rho_0 = u_0 = 1$. We note that the assumption concerning the ratio of viscosities for simple gases still holds under this change of coordinates, $16 \bar \nu  = \bar \alpha (27\Gamma+12)$. Hereafter, we drop the bar notation. To map out 
the feasible set, we solve the profile equation \eqref{prof} (with $u_{0}=1$) as an initial value problem on the interval $[0,1]$ with initial conditions $(u,e)(0)=(1,e_0)$ for various values of the integration constants $c_1$, $c_2$. We center the map about the integration constants corresponding to the fixed point, $c_1 = -1-\Gamma e_0$, $c_2 = -\frac{1}{2}-(1+\Gamma) e_0$. We check to ensure that $u$ and $e$ remain positive throughout the unit interval and that finite blowup does not occur; see Appendix \ref{appendix:feasible} for details about computational algorithms. In Figure \ref{fig:feasible}, we plot some examples of the feasible set. Note that the feasible set is unbounded (see Lemma \ref{C_unbounded}). 

We tested the following parameters to see if they lie in the feasibility set,
\begin{equation}
\begin{split}
(\Gamma,\alpha,e_0,\Delta c_1,\Delta c_2) & \in \{ 2/3,2/5,1\}\times \textrm{lin}(0.1,2,10)\times \textrm{lin}(0.001,10,30)\\ & \times \textrm{lin}(-50,50,50)\times\textrm{lin}(-50,50,50),
\end{split}
\notag
\end{equation}
where $\textrm{lin}(a,b,c)$ indicates the set containing $c$ evenly spaced points in the interval $[a,b]$, $\nu = \frac{\alpha(27\Gamma+12)}{16}$, and $c_1 =  -1-\Gamma e_0+\Delta c_1$, $c_2 = -\frac{1}{2}-(1+\Gamma) e_0+\Delta c_2$.

\subsection{Evans function computations} \label{sec:evans_computations}

To numerically compute the Evans function, we use the package STABLAB \cite{Barker_matlab}, which is well tested by this point; for example see \cite{BHLynZ1,BJNRZ,BLZ}. We provide details about numerical conditioning and algorithm choices we use in STABLAB in Appendix \ref{appendix:evans}.
\medskip

\subsection{Winding number computations}\label{s:numspec}

\begin{figure}[htbp]
 \begin{center}
$
\begin{array}{lr}
(a) \includegraphics[scale=0.39]{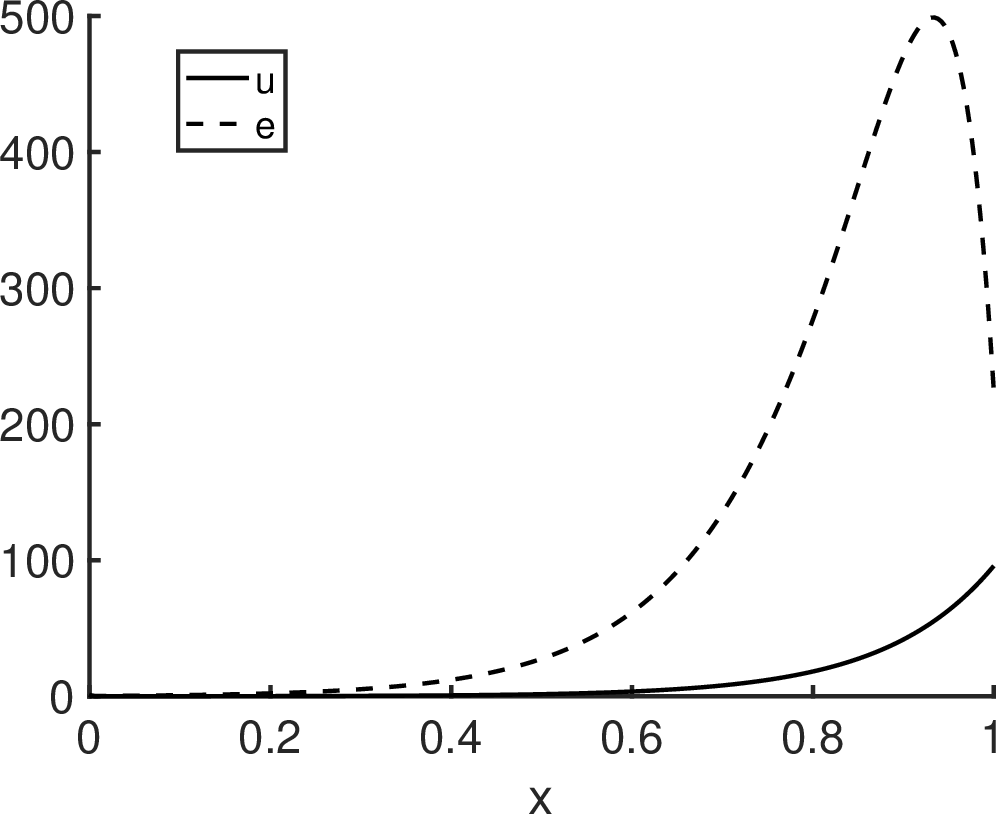} & (b) \includegraphics[scale=0.39]{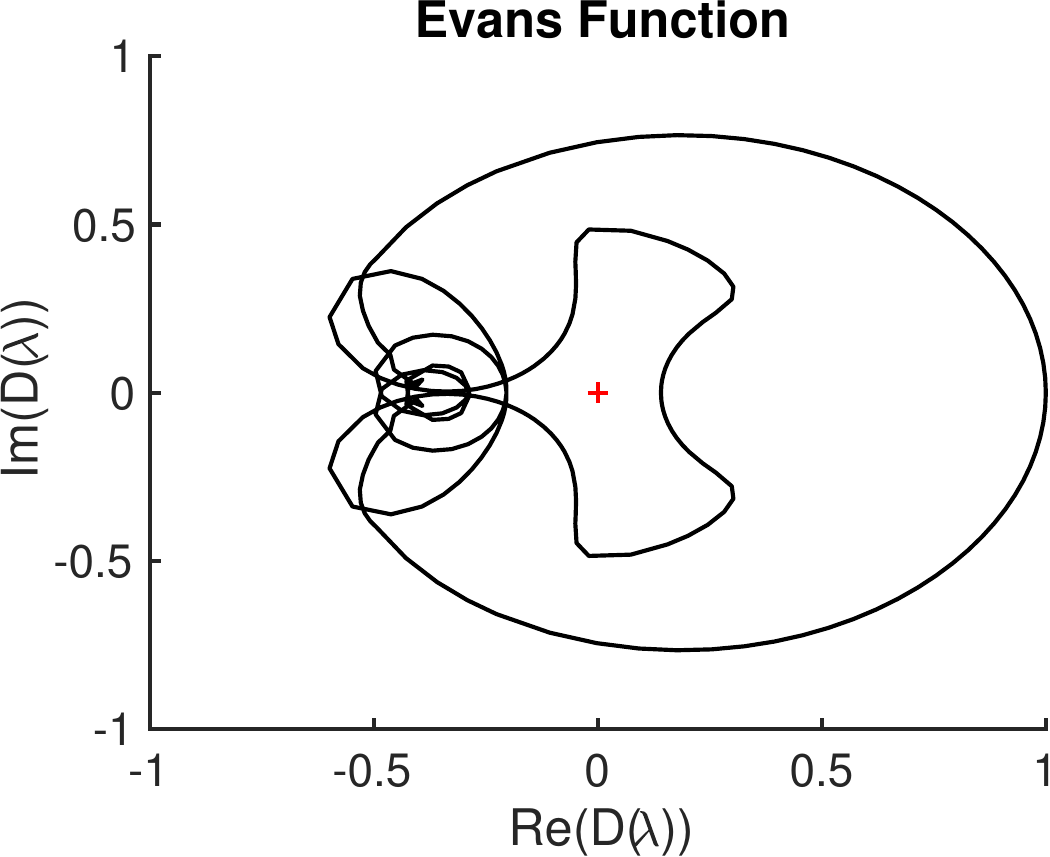}
\end{array}
$
\end{center}
\caption{For the parameters $\alpha = 0.1$, $\Gamma = 1$, $\nu = 0.2438$, $e_0 = 0.001$, $c_1 = -18.35$, and $c_2 = 0.5184$, we plot (a) the boundary layer profile, and (b) the image of $S:= \{z\in B(0,100): \Re(z)\geq 0\}$ under the Evans function. The winding number is zero indicating spectral stability of the boundary layer profile.
}
\label{fig288}
\end{figure}

To test for the existence of unstable eigenvalues, we compute the Evans function on a contour consisting of the boundary $\partial S$ of the set $S:= \{z\in B(0,100): \Re(z)\geq 0\}$. We use the functionality built into STABLAB \cite{Barker_matlab} that adaptively chooses the mesh along $\partial S$ so that the relative error between any two consecutive points on the image of $\partial S$ under the Evans function, $C_S$, varies by no more than 0.2. We then compute the winding number of $C_S$, which is the number of eigenvalues of \eqref{eval} inside $S$. In Figure \ref{fig288}, we demonstrate the profile and corresponding Evans function computation for representative parameters.

We compute the Evans function on the  contour $\partial S$ for the parameters, if they are in the feasible set, given by
\begin{equation}
\begin{split}
(\Gamma,\alpha,e_0,\Delta c_1,\Delta c_2) & \in \{ 2/3,2/5,1\}\times \textrm{lin}(0.1,2,10)\times \textrm{lin}(0.001,10,30)\\ & \times \textrm{lin}(-50,50,50)\times\textrm{lin}(-50,50,50),
\end{split}
\notag
\end{equation}
where $\textrm{lin}(a,b,c)$ indicates the set containing $c$ evenly spaced points in the interval $[a,b]$, $\nu = \frac{\alpha(27\Gamma+12)}{16}$, $c_1 =  -1-\Gamma e_0+\Delta c_1$, and $c_2 = -\frac{1}{2}-(1+\Gamma) e_0+\Delta c_2$. In all, we computed the Evans function on 670,926 contours, and in all cases found the winding number to be zero. 
These computations
took the equivalent of approximately 83.8 computation days  on a desktop with 10 duo cores.

\subsection{Computations in original coordinates}

\begin{figure}[htbp]
 \begin{center}
$
\begin{array}{lr}
(a) \includegraphics[scale=0.3]{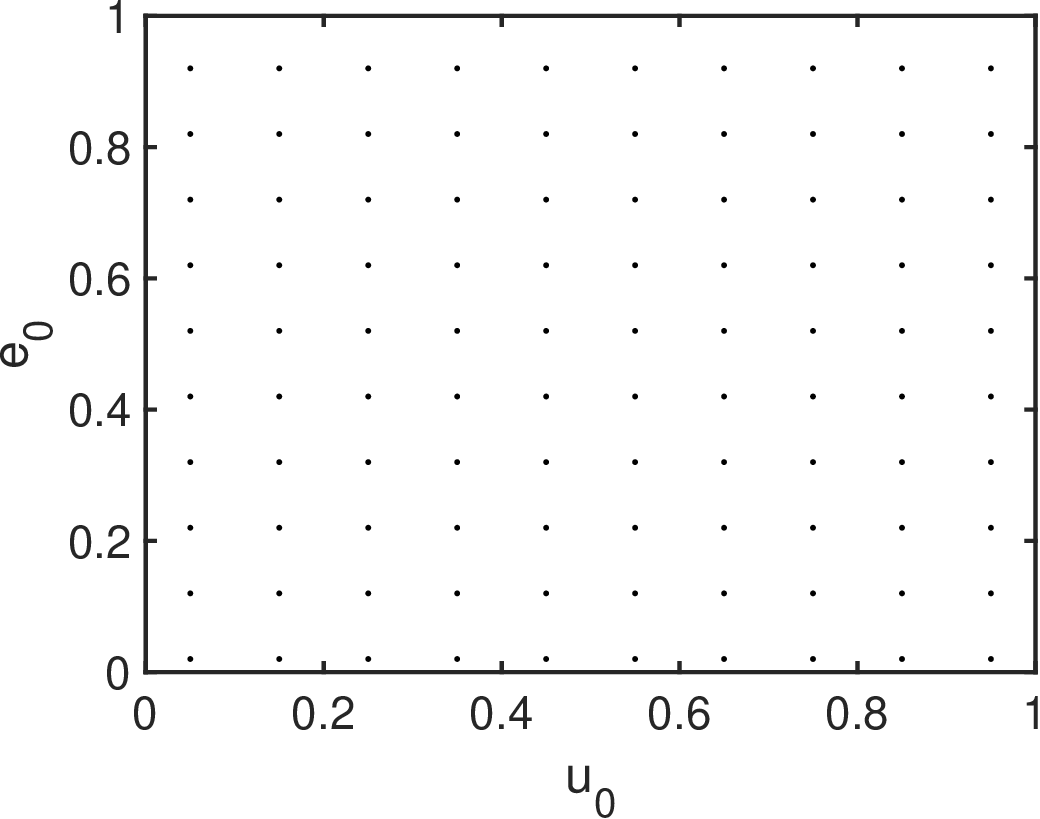} &(b) \includegraphics[scale=0.3]{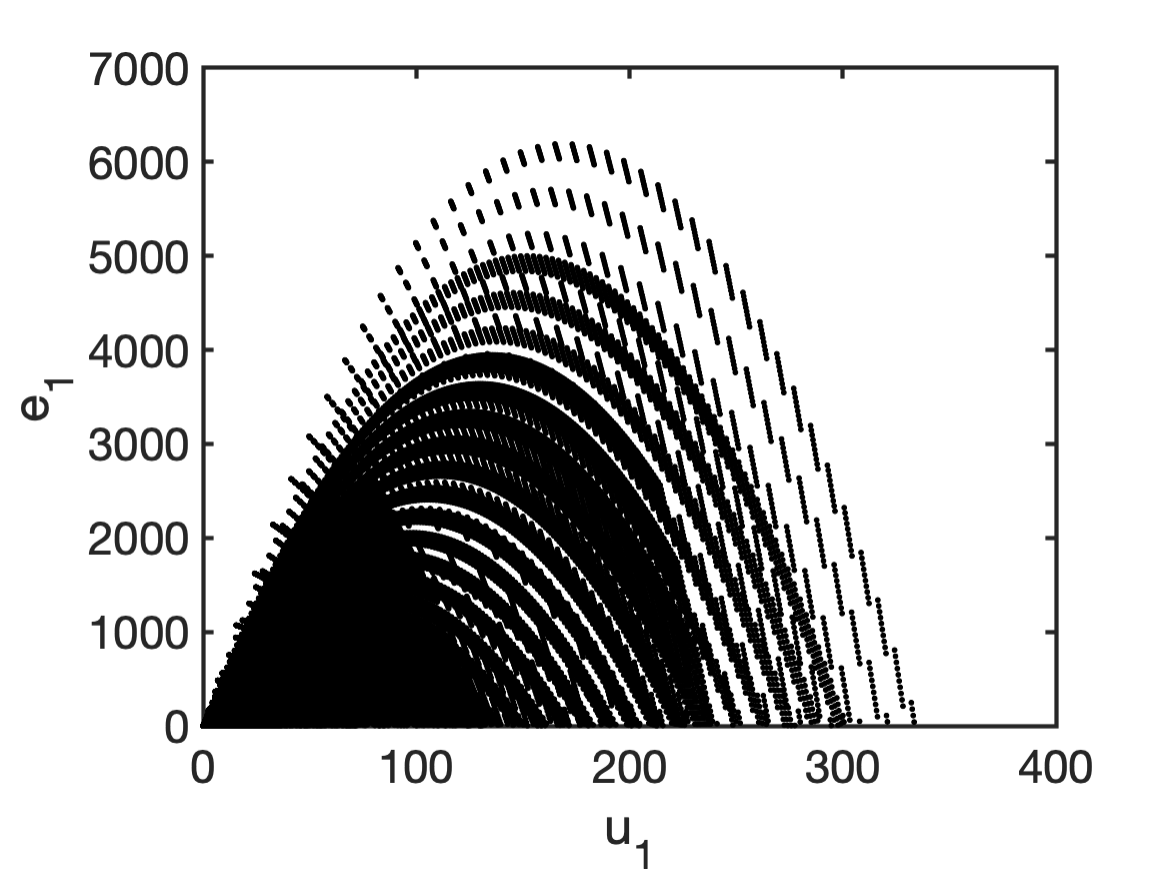}\\
(c) \includegraphics[scale=0.3]{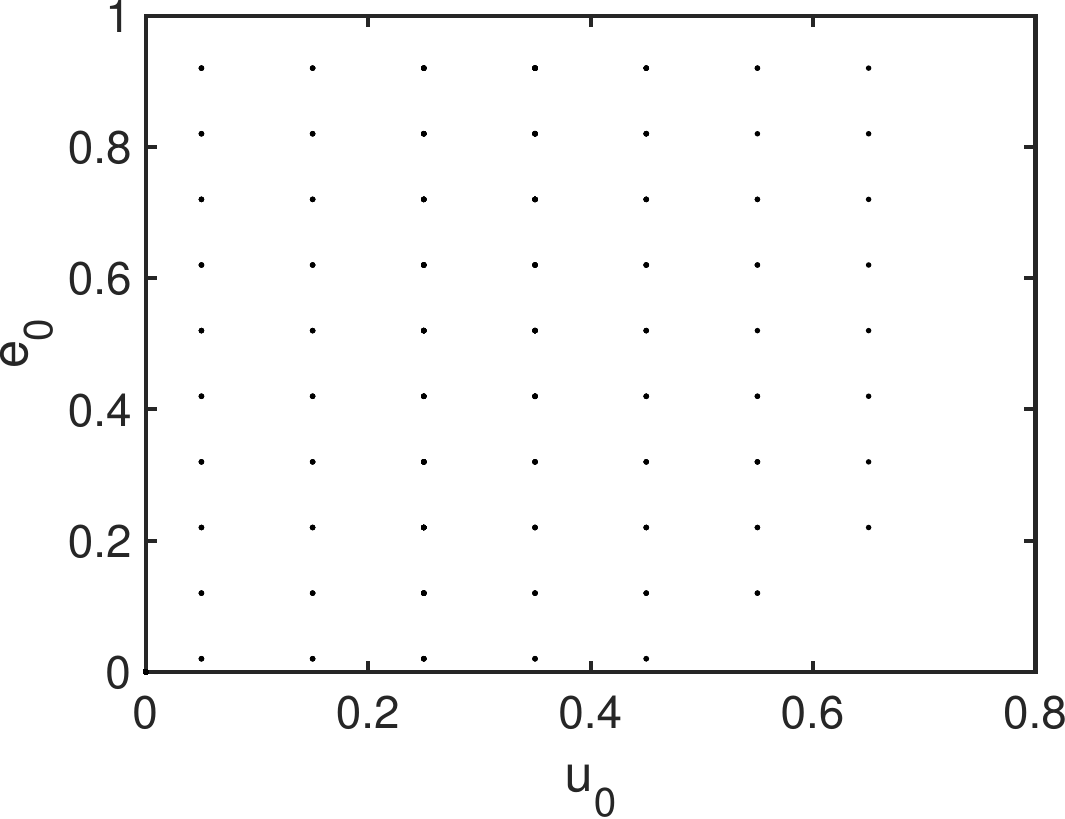} & (d) \includegraphics[scale=0.3]{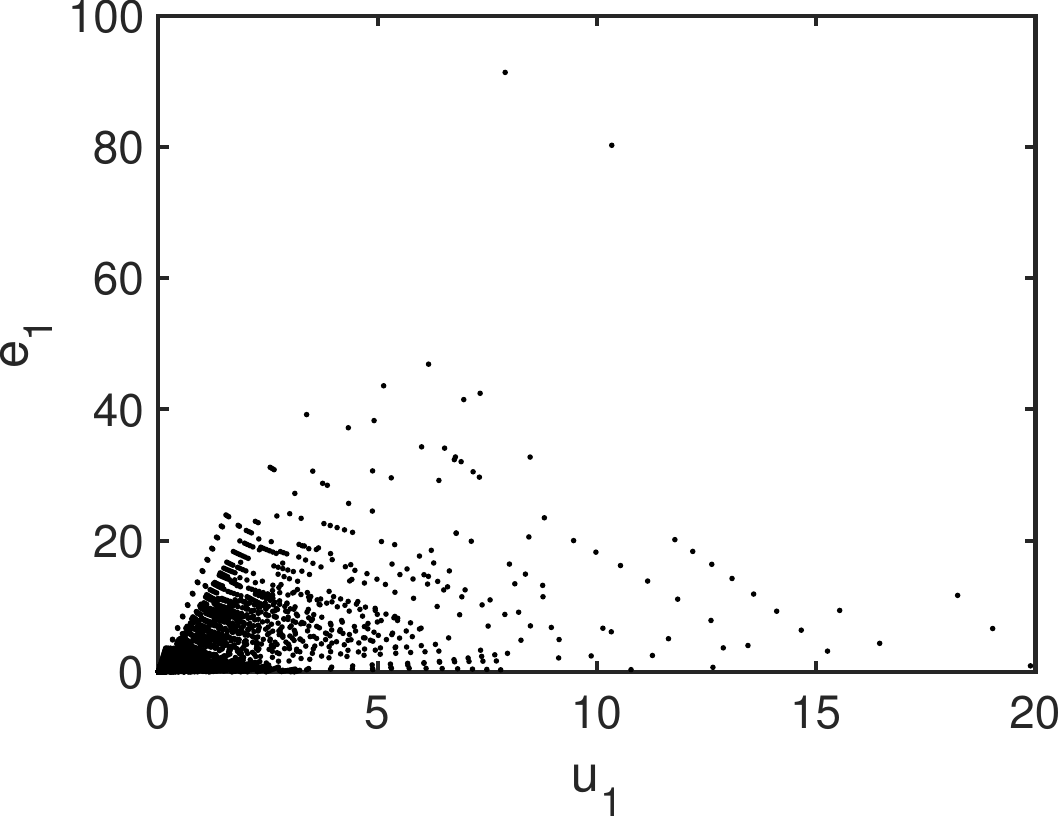}
\end{array}
$
\end{center}
\caption{(a) Plot of initial conditions $(u_0,e_0)$ considered. (b) Plot of resulting $(u_1,e_1)$ end conditions when $\Gamma = 2/3$ and $-50\leq c_1\leq 50$ (approximately 50 points) and $-50\leq c_2 \leq 50$ (approximately 200 points). In Figures (c) and (d) we plot the profile data for which we were able to compute the Evans function. }
\label{fig339}
\end{figure}

As mentioned at the beginning of Section \ref{s:num}, we use a convenient scaling for the numerics. To give an idea of the region studied in the original coordinates corresponding to the analytical results, we provide the following plots. In Figure \ref{fig339} (a), we plot the coordinates for the initial data, and in Figure \ref{fig339} (b) we plot the coordinates for the final data, for the profiles for which existence is shown. In Figures \ref{fig339} (c) and (d), we plot the same initial and final profile data for which we were able to compute the Evans function and show stability. 
We note that computing the Evans function, though straightforward on much of the feasibility region,
is very challenging in some parts due to stiffness of the associated ODEs. 

This stiffness can come from rapid transition of the profile solution associated with differences in initial and final data that are $\gg 1$; see Figures \ref{fig339} (a)-(b) for examples of such data. It can also come from
points near the feasibility boundary, where the profile nearly blows up to infinity or down to zero 
on the interior of the computational interval; as one can see from the linearized profile equations \eqref{eqD0},
either of these events leads to blowup of coefficients and associated stiffness of this system of ODE.

\subsection{Global uniqueness/stability index}\label{s:global}

\begin{figure}[htbp]
 \begin{center}
$
\begin{array}{lr}
(a) \includegraphics[scale=0.3]{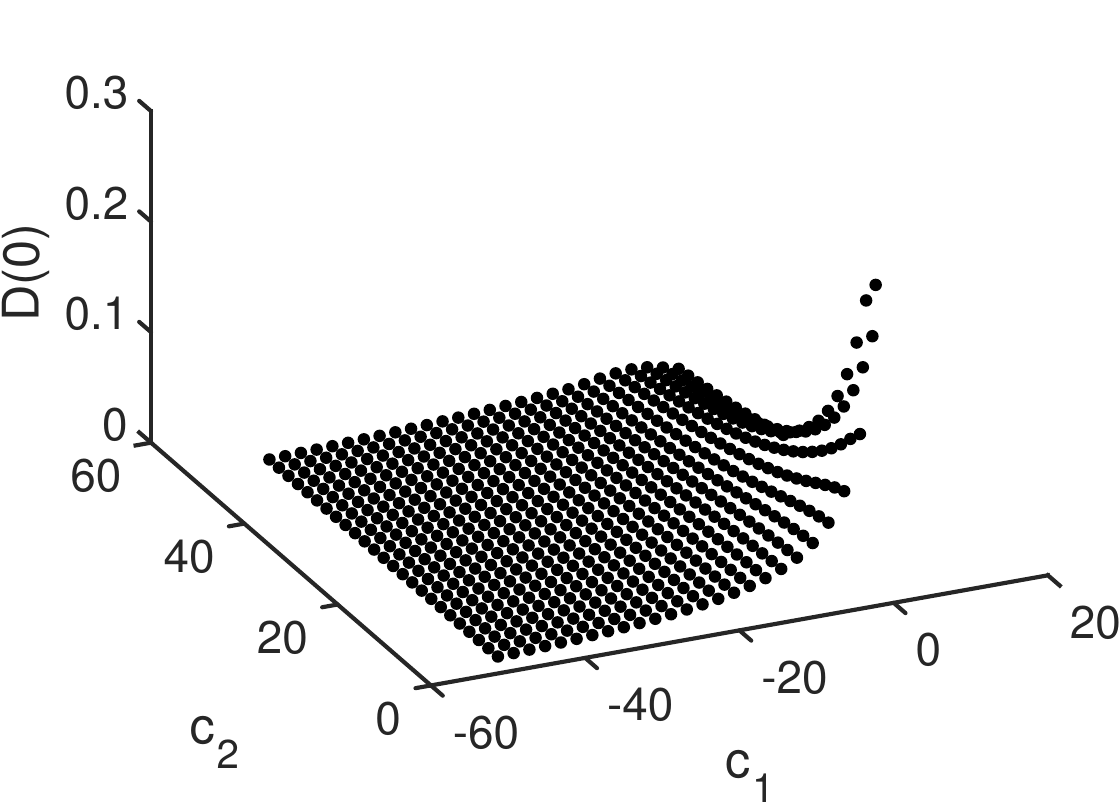} & (b) \includegraphics[scale=0.3]{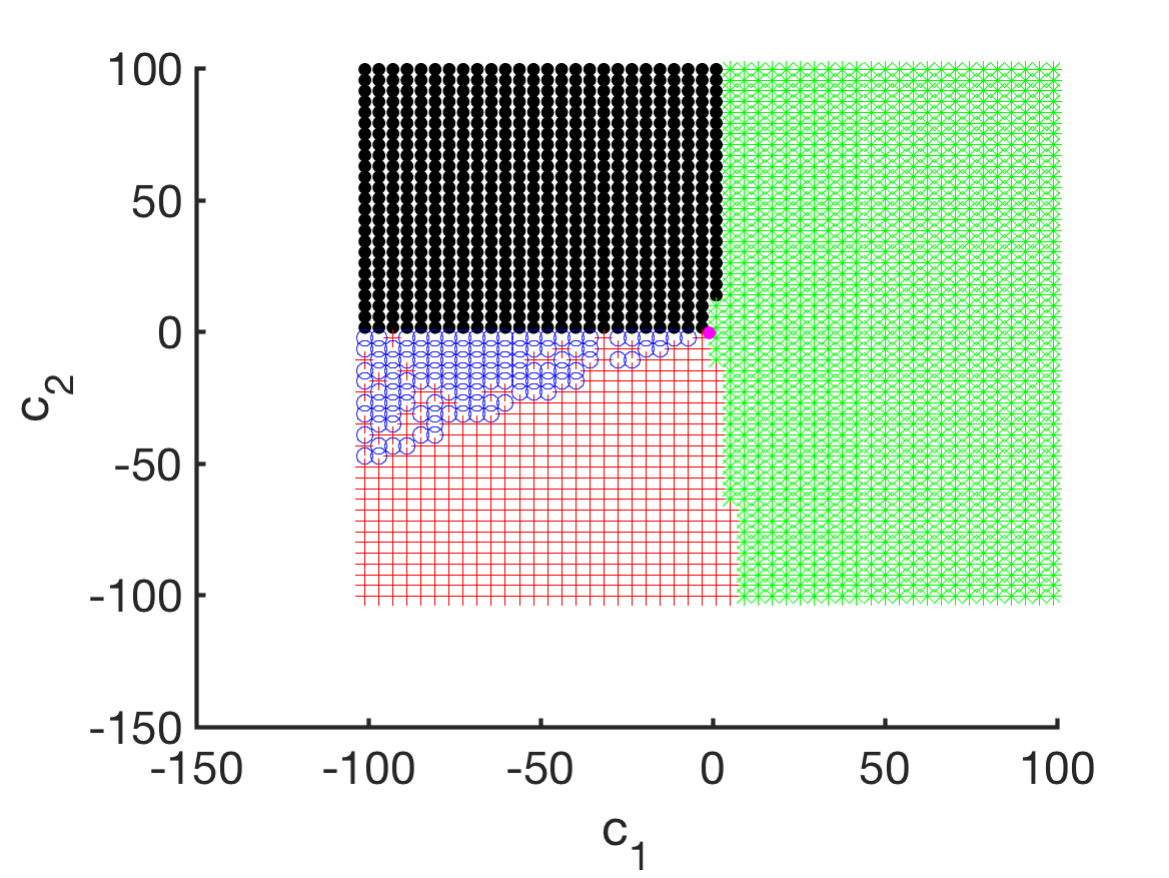}
\end{array}
$
\end{center}
\caption{In these figures, $\alpha = 0.7\bar 3$, $\Gamma = 2/3$, $\nu = 1.375$, and $e_0 = 0.001$. (a) Plot of $D(0)$ against $c_1$ and $c_2$. (b) Plot of the feasible set corresponding to Figure (a) with black dots, the set where $u$ goes negative on $[0,1]$ with blue circles, the set where $e$ goes negative on $[0,1]$ with green stars, and the set where there is finite time blowup on $[0,1]$ with red $+$ signs. 
}
\label{fig283}
\end{figure}

For the parameters in the feasible set described in Section  \ref{s:numfeas}, we computed the Evans function at the origin, $D(0)$. We found that the smallest value of $D(0)$ for the computed parameters was 4.19e-4,
far greater than the absolute error tolerance 1e-10 of the computation method, as described further in
Appendix \ref{appendix:non_unique}.
Thus, $D(0)$ appears not to vanish on the feasible set; 
See Figure \ref{fig283} for a demonstration of how $D(0)$ varies as $c_1$ and $c_2$ vary in the feasible set.
By Proposition \ref{uniqueprop}, uniqueness of profiles is equivalent to
nonvanishing of the Jacobian determinant $\gamma$ associated with the profile ODE,
while Lemma \ref{ZSlem} shows that $\gamma$ is a real positive multiple of $D(0)$.
Thus, the numerically observed nonvanishing of $D(0)$ indicates
global uniqueness on the feasible set, as announced in the introduction.

Moreover, nonvanishing of $D(0)$ implies further that the sign of the
stability index $\mu$, defined in \eqref{mudef} as the sign of a nonvanishing real multiple of $D(0)$.
is constant. 
By numerical evaluation at a single choice of profile-- alternatively,
by continuation to the constant-profile limit, for which the limiting sign of $D(\lambda)$ at positive
real infinity may be computed to be $\mu=+1$-- we find that $\mu\equiv +1$ on the feasible set, consistent with
(though not implying) spectral stability.\footnote{ 
That $\mu=+1$ for constant profiles may be deduced also by the fact \cite[Proposition 3.2]{paper2} 
that constant steady solutions of general systems with convex entropy are stable.}

\section{A numerical counterexample }\label{s:ceg}

We now consider equations \eqref{nsint} subject to the equation of state 
$\bar e(\tau,S) = \frac{e^S}{\tau} + S+\frac{\tau^2}{2}$
considered in \cite{BFZ}, where $\tau$ corresponds to specific volume and $S$ corresponds to entropy. Specific density is given by $\rho = \frac{1}{\tau}$, and $T = \bar e_S = \frac{e^S}{\tau} + 1$, so 
$e^S = \frac{T-1}{\rho}$, or
\be\label{localent}
S =\hat S(\rho, T)=\ln \left( \frac{T-1}{\rho} \right).
\ee
From this, we obtain 
\be\label{localpress}
\begin{split}
p &=\hat p(\rho,T)= -\bar e_{\tau} = \rho(T-1)-\frac{1}{\rho},
\\ e&=\hat e(\rho,T)=T-1 + \ln \left( \frac{T-1}{\rho} \right) + \frac{1}{2}\rho^2,
\end{split}
\ee
closing the system, together with the energy relation $E=e+\frac{1}{2}u^2$, in terms of variables
$(\rho, u, T)$, $T>1$.
Alternatively, inverting the relation $e=\hat e(T,\rho)$ using $\hat e_T>0$ for $T>1$, 
we may consider it as implicitly determining a system in the usual variables $(\rho, u,e)$,
with $e>0$.
This is the system referred to as the \textit{local model} in \cite{BFZ}.
Notably, the function $\eta:=-\rho \hat S(\rho, T)$, with $\hat S$ as in \eqref{localent}
considered as a function of the conservative variables $(\rho, \rho u, E)$ 
{\it is a convex entropy for system \eqref{nsint}} in the sense of 
\cite{La,KSh}; see \cite{BFZ}.

In \cite{BFZ} it was shown that the local model considered on the whole line 
has unstable shock waves for parameters for which the inviscid system has stable waves. 
Here, we demonstrate that the local model considered on a finite interval
has parameters for which uniqueness of solutions fails, and also other, nearby parameters 
for which a Hopf-bifurcation occurs.

These results are guided by the general principles of \cite{Z1}
relating spectra of standing shocks on the whole line to spectra of pieces thereof, 
considered as solutions on a truncated domain.
See \cite[\S 3.2]{paper2} for further discussion in the specific case of a finite interval.
The first relevant principle is that spectra on the interval are, for $\Re \lambda \geq 0$ and $\lambda \neq 0$ given in
the limit as interval length goes to infinity- equivalently, as viscosity goes to zero- 
by the direct sum of spectra on the whole line together with spectra of constant boundary layers on the half-line
with data corresponding to that on the left (resp. right) endpoint of the interval.  This implies that strict
instability on the whole line implies strict instability on the interval 
with associated stability transition as amplitude is increased from a (presumably stable; see
\cite[Proposition 3.2]{paper2}) constant steady solution to an 
unstable one.

The second principle is that in the same large interval length/small viscosity standing-shock limit, 
the stability index does not vanish (\cite[Prop. 3.3]{paper2}), or equivalently $D(0)\neq 0$.
Thus, if a homotopy is taken from 
stable constant solutions to unstable standing shock solutions, entirely within the class of standing shocks with
sufficiently large interval/small viscosity, then the associated stability transition cannot correspond to a simple
crossing of an eigenvalue through the origin $\lambda=0$, as $D(0)\neq 0$, and must therefore involve the crossing of
one or more pairs of complex conjugate roots, i.e., a Hopf-type scenario.

On the other hand, the first cited principle implies that two of these roots must be near the pair of roots at the origin
of the whole-line shock as it undergoes transition to instability: one ``translational'' eigenvalue fixed at $\lambda=0$
and the crossing eigenvalue corresponding to instability.  Thus, we have the picture of a Hopf bifurcation with 
very nearby roots, i.e., with associated time-period going to infinity, a quite delicate scenario.
This makes numerical verification somewhat sensitive; however, it also aids us in finding a more standard bifurcation in
the form of a single crossing eigenvalue through $\lambda=0$, as we are able to find by playing with the left and right 
boundaries of the interval for a given, sufficiently large-amplitude standing shock on the whole line.

\subsection{Nonuniqueness}\label{s:nonuniqueness}

\textbf{Abstract bifurcation result.}
We first demonstrate (numerically) a bifurcation implying nonuniqueness.
Namely, considering the restriction to a finite interval $[x_L,x_R]$,
of an appropriate standing-shock solution of the local model on the whole line (described in detail below),
we show that $D(0)$ changes sign as $x_L$ and $x_R$ vary; see Figure \ref{eq:bfz_nonuniqueness} (a)-(c). 
Defining by $c_*(x_L, x_R)$ the value of $c$ corresponding to the shock profile on $[x_L,x_R]$, define the map
$$
\Phi(c; x_L,x_R):= \psi(c_*(x_L,x_R) + c; x_L, x_R)- \psi(c_*(x_L,x_R); x_L, x_R),
$$
where $\psi(x;x_L,x_R)$ is the solution map $\psi$ associated with the interval $[x_L,x_R]$.
Then $\Phi(0;x_L,x_R)\equiv 0$, reflecting the fact that the shock profile restricted to $[x_L,x_R]$
solves its own data.  Existence of additional roots $c\neq 0$ for some $x_L$, $x_R$ implies nonuniqueness for the
same data.
Nonuniqueness may be detected, therefore, using the following abstract bifurcation result, in the spirit 
of Proposition \ref{uniqueprop} and \cite[Lemma 3.10]{BFZ}.

\bpr\label{bifprop}
Let $\Phi(c; p):\R^m\times \R$ satisfy $\Phi(0;p)\equiv 0$. If $\gamma:=\det (d\Phi(0;p))$ 
changes sign as $p$ crosses a particular bifurcation value $p=p_*$, then $\Phi(\cdot;p)$ has a nontrivial
root $c\neq 0$ for $p$ arbitrarily close to $p_*$.
\epr

\begin{proof}
Arguing by contradiction, suppose that $c=0$ is the unique root of $\Phi(c;p)=0$ for $p$ in a neighborhood
of $p_*$.  Thus, $\Phi$ does not vanish on the boundary of a small ball $B(0,r)$, hence the topological degree
of $\Phi(\cdot;p)$ is independent of $p$.  However, at $p$ for which $\det (d\Phi(0;p))>0$, the degree is by the
assumed uniqueness of roots equal to $+1$, while at points $p$ for which $\det (d\Phi(0;p))<0$, the degree is $-1$,
a contradiction.
\end{proof}

 \medskip

\textbf{Numerical investigation: methods.}
To demonstrate non-uniqueness numerically, we first solve approximately
for the profile corresponding to the whole-line shock. 
We then take the piece of that solution on $[x_L,x_R]$ as the profile for the finite boundary problem posed on the same interval. 
Next, following the intuition described at the end of the previous subsection, we seek nonuniqueness
by appropriately varying parameters of this ``truncated shock'' and nearby steady solutions.

The computations showing non-uniqueness are relatively difficult. In the following discussion, $S_- := \lim_{x\to -\infty} S(x)$, is the left end state value of entropy in the whole-line shock wave solution of the  local model.
 To solve for the profile, we fix the parameters $\alpha = \kappa = 1$ and take $S_- = 1$. From the Rankine-Hugoniot conditions, we obtain the other parameters. We then use a boundary value solver to obtain the whole-line viscous shock solution. 
 
 Next,  we use continuation with 30 evenly spaced steps in $S_-$ to obtain the solution at $S_- = -5$. That is, we change the parameter $S_-$ by a small amount and solve for the other parameters given by the Rankine-Hugoniot conditions, then use the profile solution corresponding to the previous value of $S_-$ as an initial guess in the boundary value solver  to solve for the profile for the new parameters. In solving for the whole-line profile, we use STABLAB which adaptively increases the spatial domain $[-L,L]$, $L \gg 1$, until the profile converges to the fixed-point end states, corresponding to the shock at $x = \pm \infty$, to within requested tolerance. 
 
To compute the Evans function, we use the same procedure as described in Section \ref{sec:evans_computations}, except that we evaluate the Wronskian to obtain the Evans function at $x = 0$ instead of $(x_L+x_R)/2$, and we use
 ``pseudo-Lagrangian coordinates'' as described in \cite{BHLynZ3} to reduce winding in our winding number studies
 without changing the zeros of the Evans function. For algorithm details, see Appendix \ref{appendix:non_unique}.

\textbf{Nonuniqueness.}
 To demonstrate abstract non-uniqueness of profile solutions, 
we take a piece of the whole-line shock for an unstable wave in the local model, 
and truncate it to a finite interval as described just above. 
 By varying the boundary on the left of this finite interval, we are able to observe a change of sign of the Evans function evaluated at the origin, $D(0)$, indicating by Proposition \ref{bifprop} that non-uniqueness of solutions occurs. 

 \medskip

\begin{figure}[htbp]
 \begin{center}
$
\begin{array}{cc}
\includegraphics[scale=0.37]{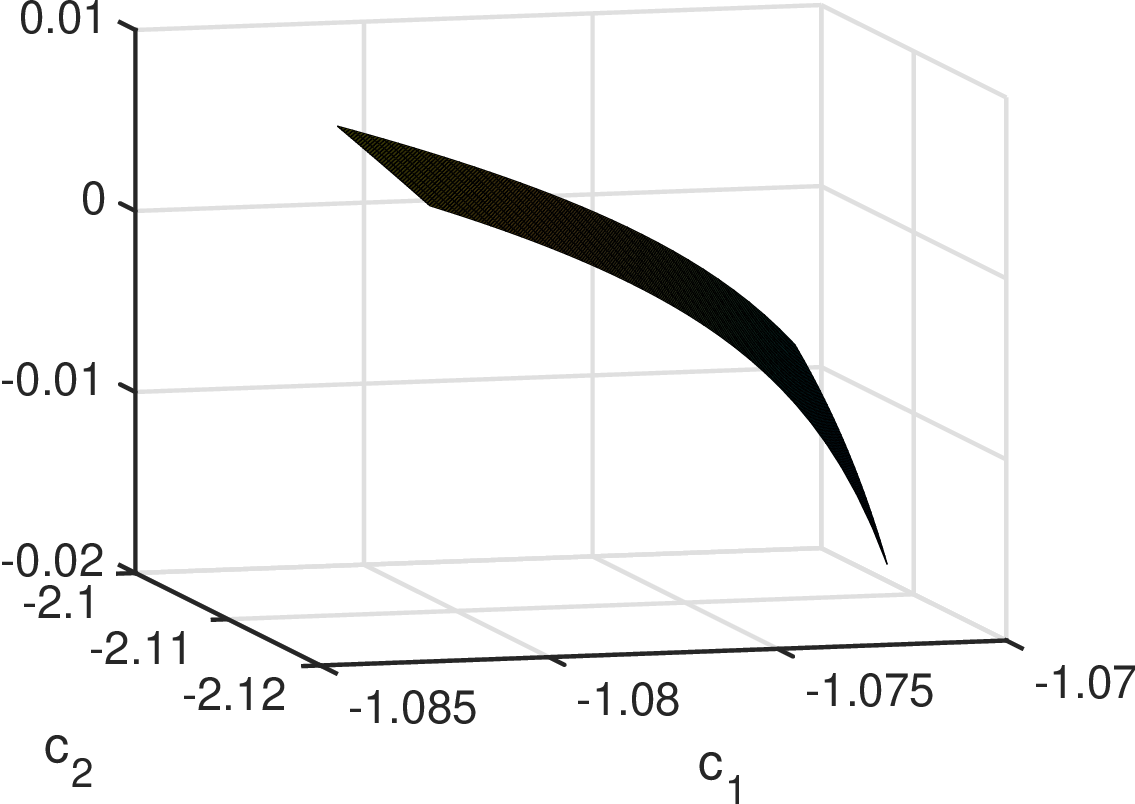} &\includegraphics[scale=0.37]{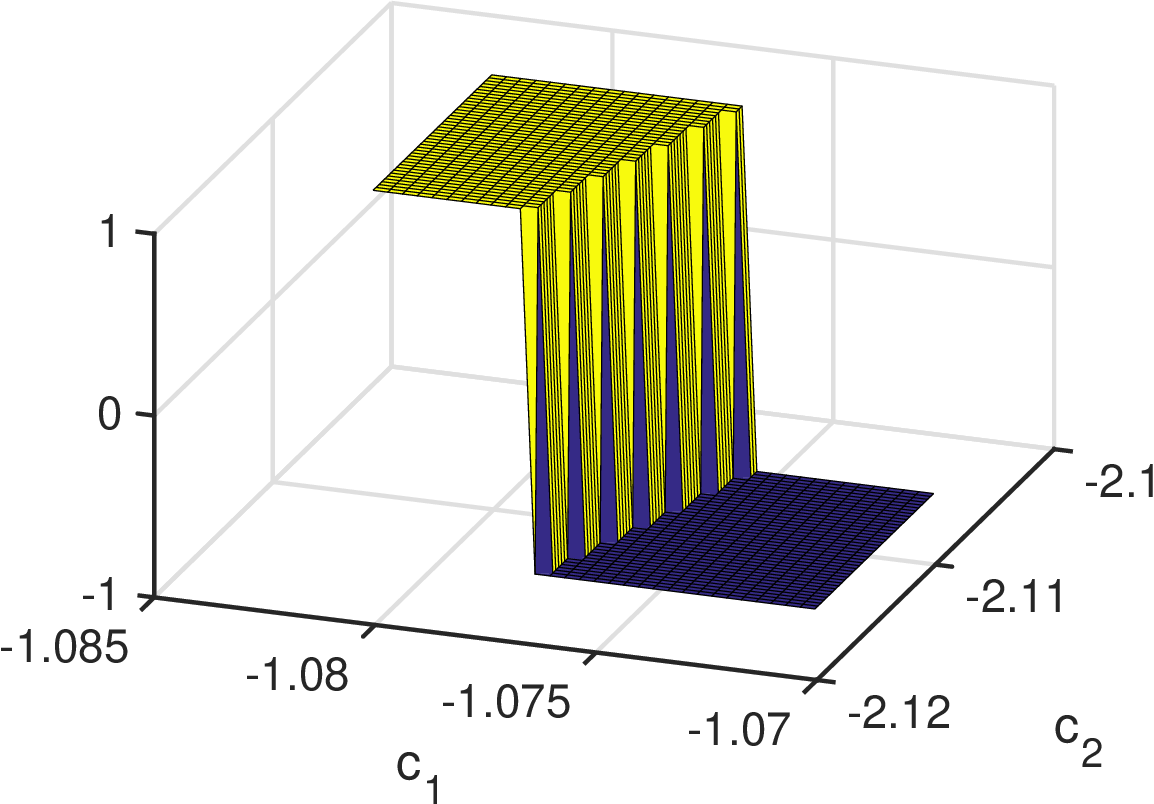}\\ 
(a)&(b)\\
\includegraphics[scale=0.37]{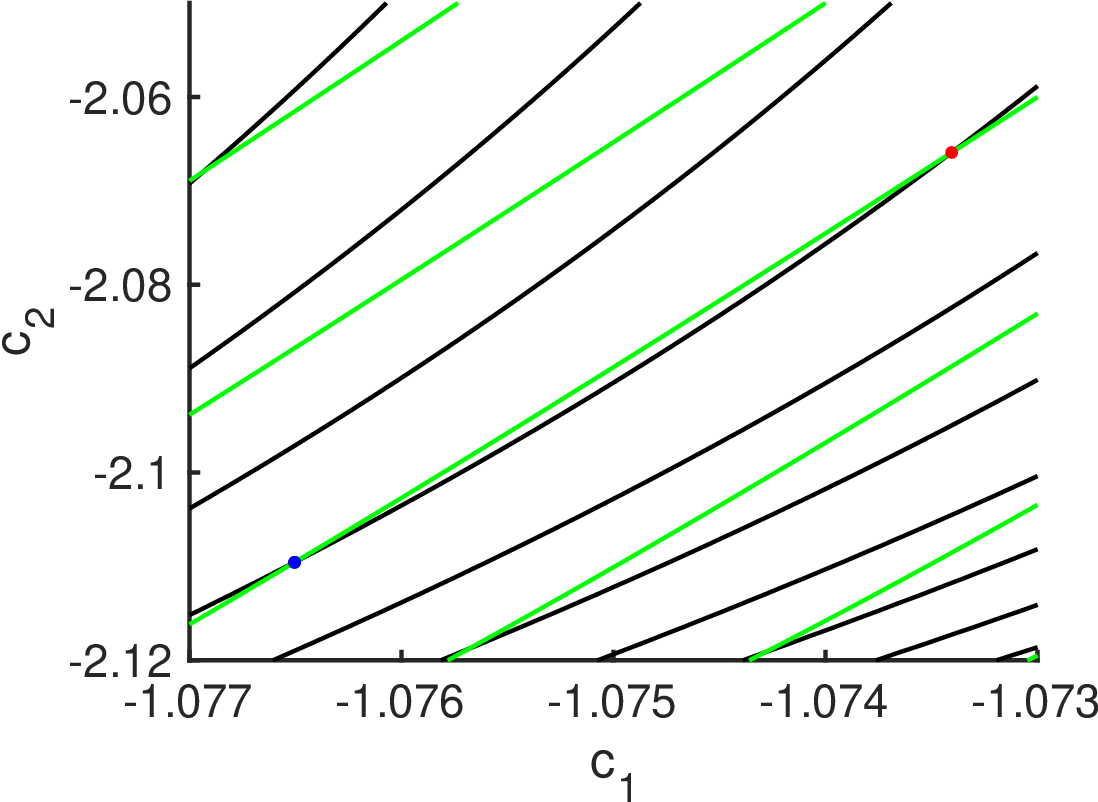} & \includegraphics[scale=0.37]{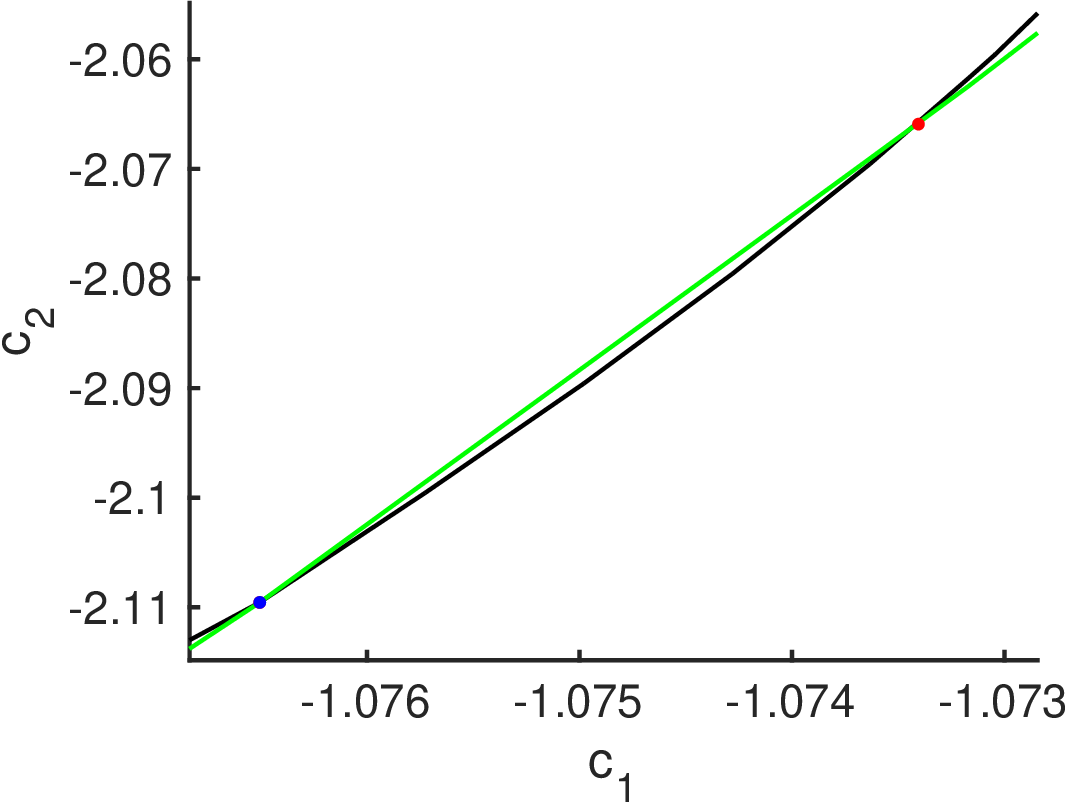}\\
(c)&(d)
\end{array}
$
\end{center}
\caption{Figures (a)-(b) demonstrate that $D(0)$ changes sign as $c_1$ and $c_2$ vary. Figures (c)-(d) indicate that there are distinct profiles that solve the same data since there are nullclines of $M_1$ and $M_2$ that intersect twice. (a) Plot of $D(0)$ against $c_1$ and $c_2$. (b) Plot of $\textrm{sign}(D(0))$ against $c_1$ and $c_2$. (c) Plot of the nullclines of $M_1$ and $M_2$. Dots indicate intersections of the nullclines. (d) Plot of only the two intersecting nullclines seen in (c).}
\label{eq:bfz_multiple}
\end{figure}

\textbf{Multiple solutions.}
%
 To find explicit nonunique profiles satisfying the same data, we investigate further, considering
 not only translated pieces of the unstable whole line shock, but also other steady solutions nearby.
 Namely, fixing the interval to be $[x_L,x_R] = [-33.17,2.9]$, 
 we compute the Evans function at the origin for profiles with varying $c_1$ and $c_2$ to find regions in $c_1$ and $c_2$ for which $D(0)$ has opposite sign; see Figures \ref{eq:bfz_multiple}(a)-(b).
 As clearly evident in Figure \ref{eq:bfz_multiple}(b), the results are consistent with a fold bifurcation of the
 map $\psi$, with orientation $\sgn \det d\psi$ changing across a smooth curve in $(c_1,c_2)$,
 in which case we may expect an open set of distinct parameter pairs $(\hat c_1,\hat c_2)$ 
 and $(\tilde c_1,\tilde c_2)$ corresponding to different profiles solving the same data, lying in regions 
 for which $D(0)$ has opposite sign. 

 A nice way to numerically find 
 parameter pairs corresponding to two distinct profiles solving the same data is to look at 
 nullclines of the mappings $M_1(c_1,c_2) := u_L(c_1,c_2)-u_L^*(c_1^*,c_2^*)$ and $M_2(c_1,c_2) := T_L(c_1,c_2)-T_L^*(c_1^*,c_2^*)$. Here $(c_1^*,c_2^*)$ are fixed constants of integration that correspond to the whole line shock, which constants of integration we find by solving for them in the Rankine-Hugoniot equation. The other terms used in defining $M_1$ and $M_2$, that is $u_L(c_1,c_2)$ and $T_L(c_1,c_2)$, are the components of the profiles evaluated at $x = x_L$. We note that these profiles have the same data at $x = x_R$ as the the profile corresponding to $(c_1^*,c_2^*)$. In particular, $u_R(c_1,c_2) = u_R(c_1^*,c_2^*)$ and $T_R(c_1,c_2) = T_R(c_1^*,c_2^*)$. 
 The level sets of $M_1$ and $M_2$ intersect in two locations, which we name $(\hat c_1,\hat c_2)$ and $(\tilde c_1,\tilde c_2)$, along the same curves, indicating that these constants of integration correspond to two distinct profiles solving the same data; see Figure \ref{eq:bfz_multiple} (c)-(d). 

\begin{figure}[htbp]
 \begin{center}
\begin{equation*}
\begin{array}{c}
(a) \includegraphics[scale=0.4]{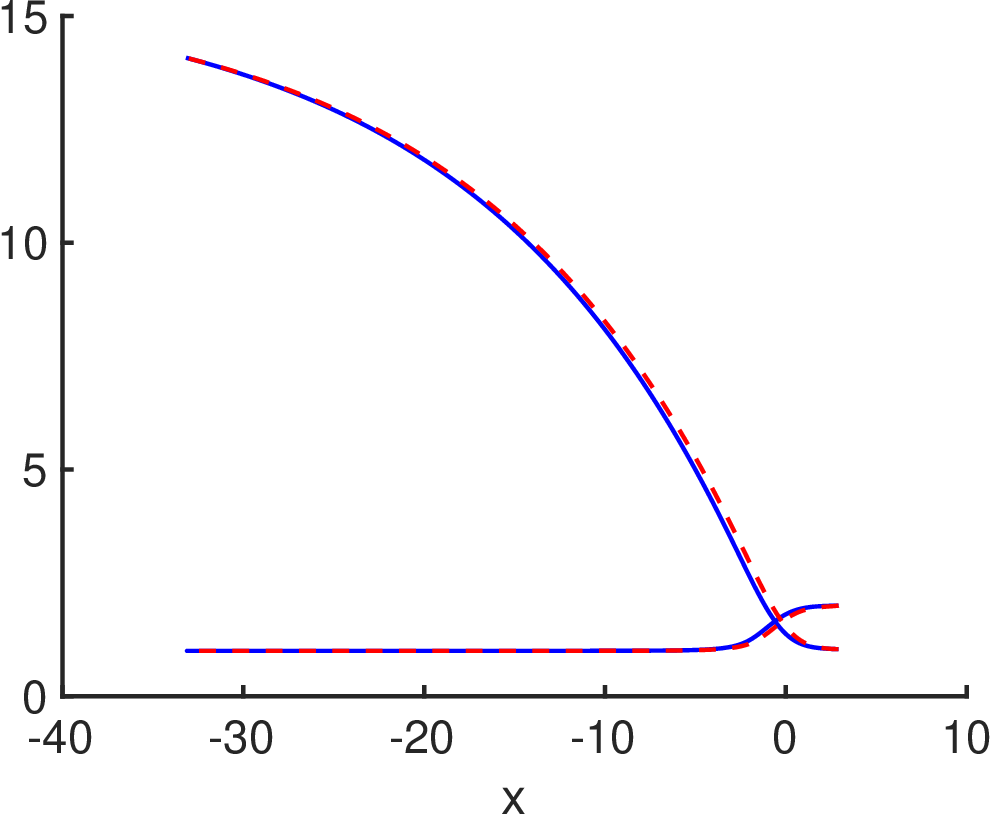} \\ (b) \includegraphics[scale=0.4]{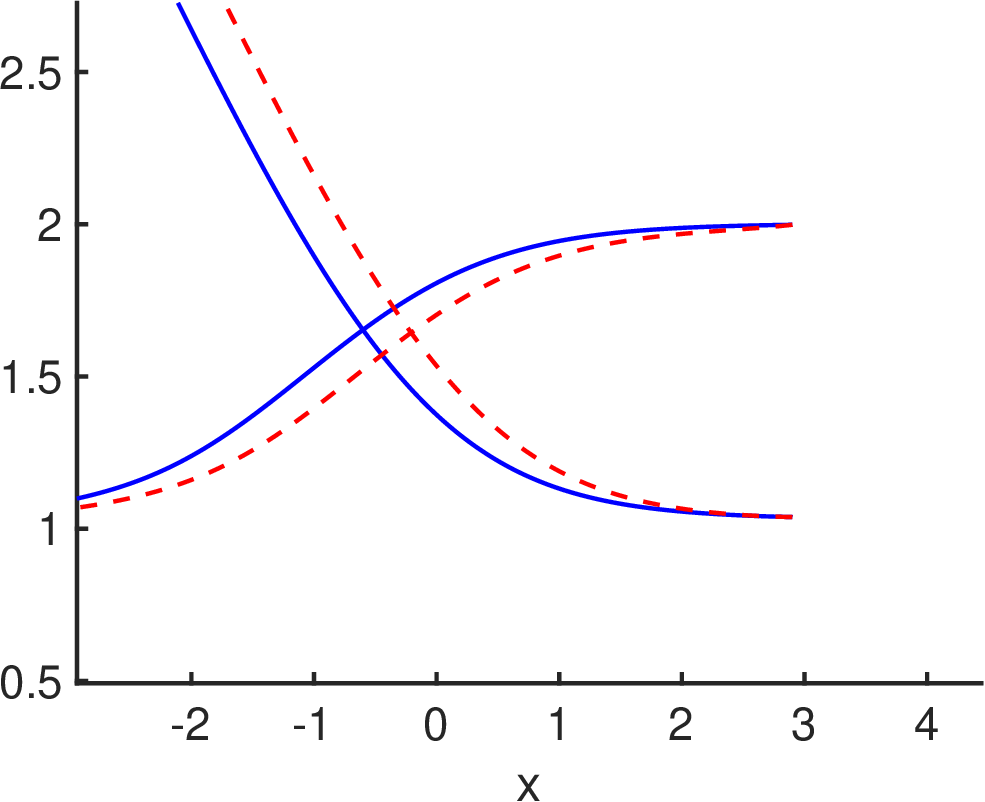}\\
\end{array}
\end{equation*}
\end{center}
\caption{(a) Plot of the two profiles, solving the same data, against $x$. The solid blue curves and dashed red curves correspond to the profiles with $c_1$ and $c_2$ values plotted as dots with the same colors in Figure \ref{eq:bfz_multiple}(c)-(d). (b) Zoomed-in picture of (a) near $x = x_R$.}
\label{eq:bfz_multiple2}
\end{figure}

\begin{figure}[htbp]
 \begin{center}
\begin{equation*}
\begin{array}{lcr}
(a) \includegraphics[scale=0.24]{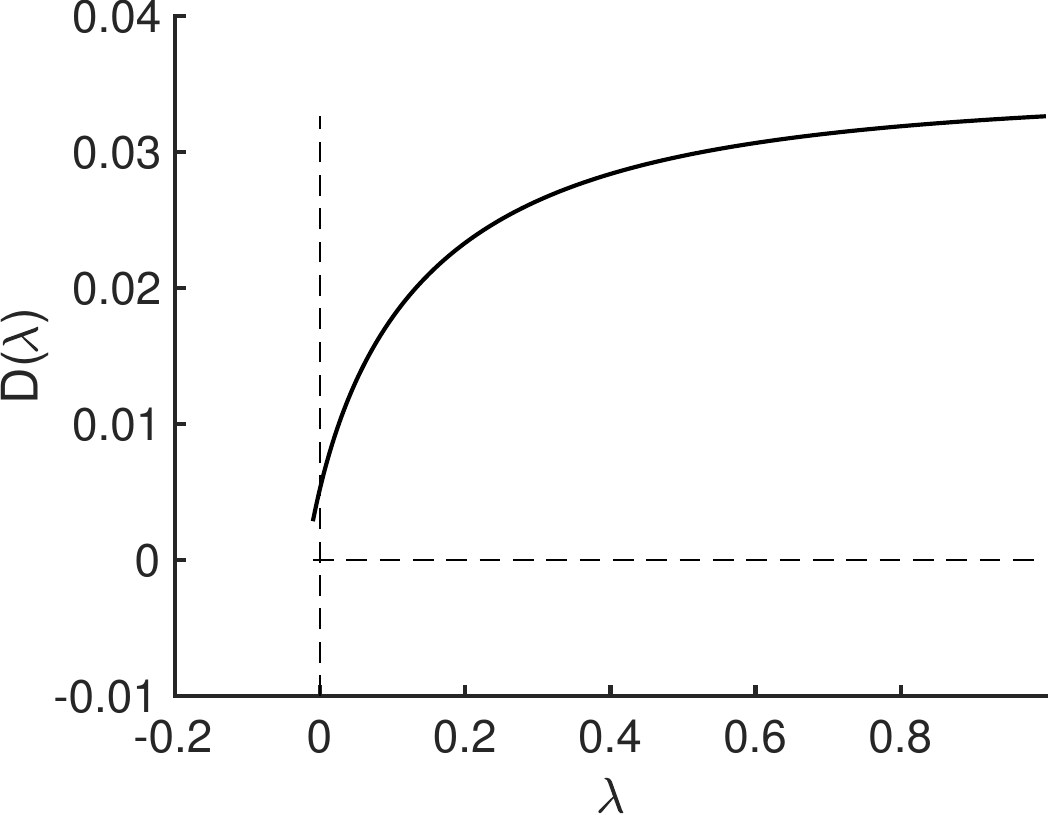} & (b) \includegraphics[scale=0.24]{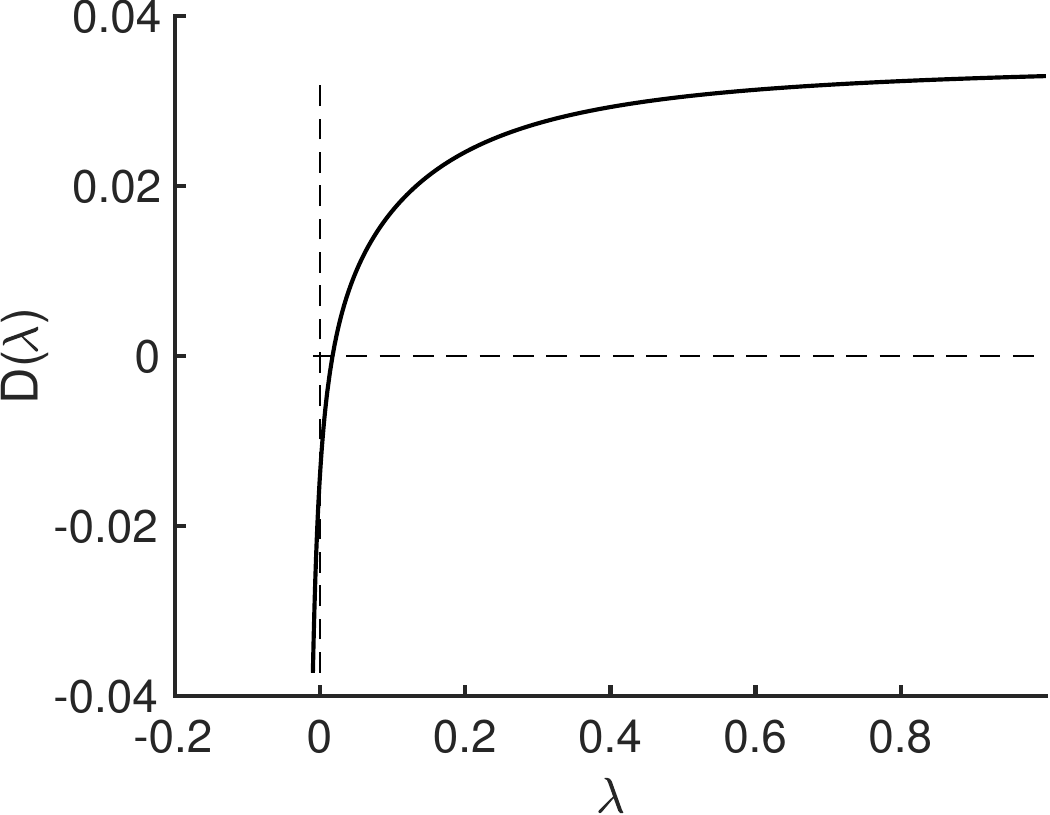} (c)  \includegraphics[scale=0.24]{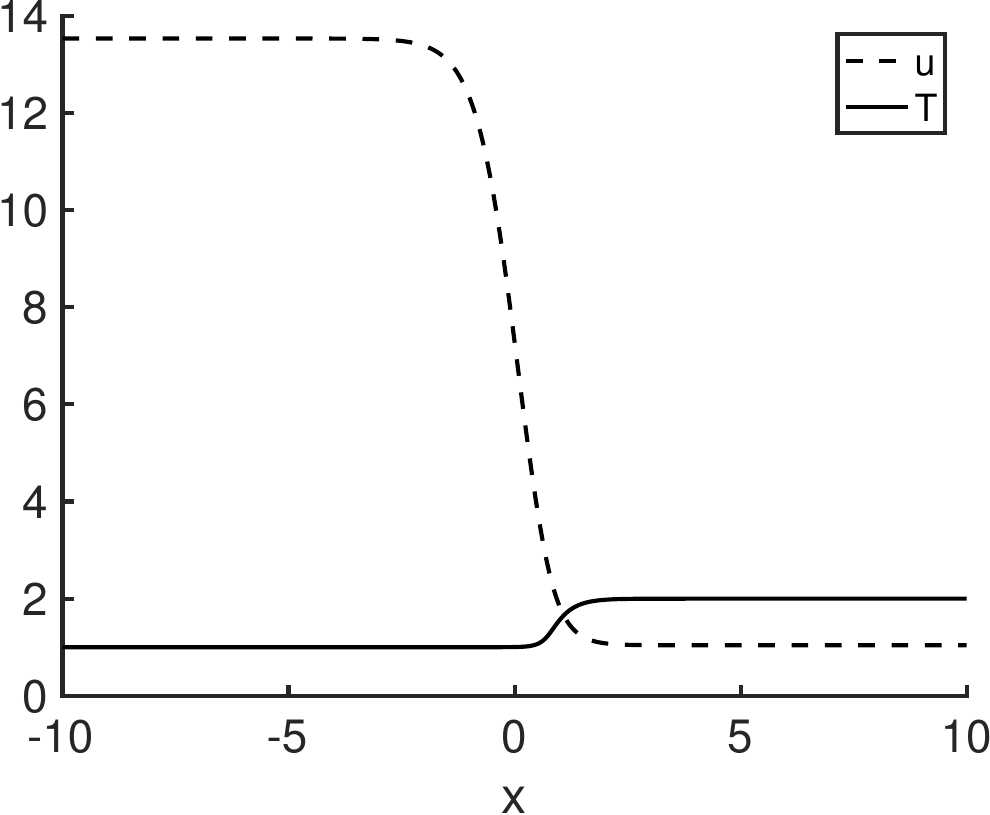}  \\
(d) \includegraphics[scale=0.24]{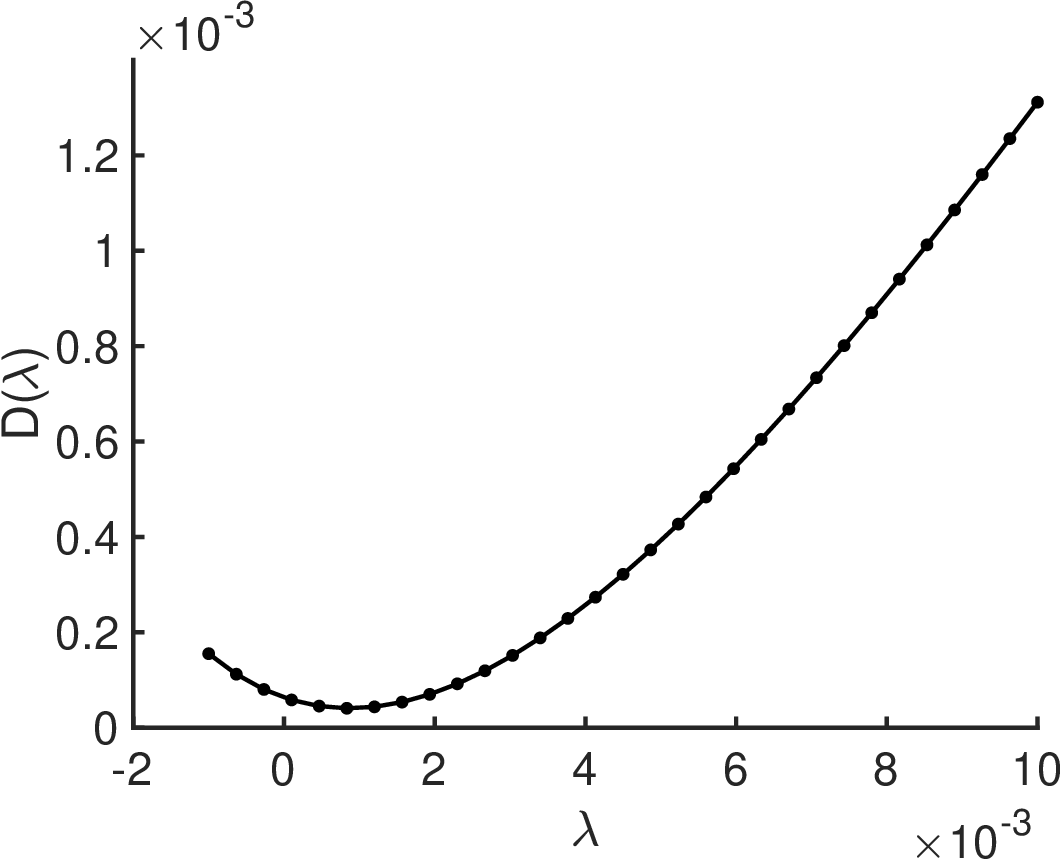} & (e) \includegraphics[scale=0.24]{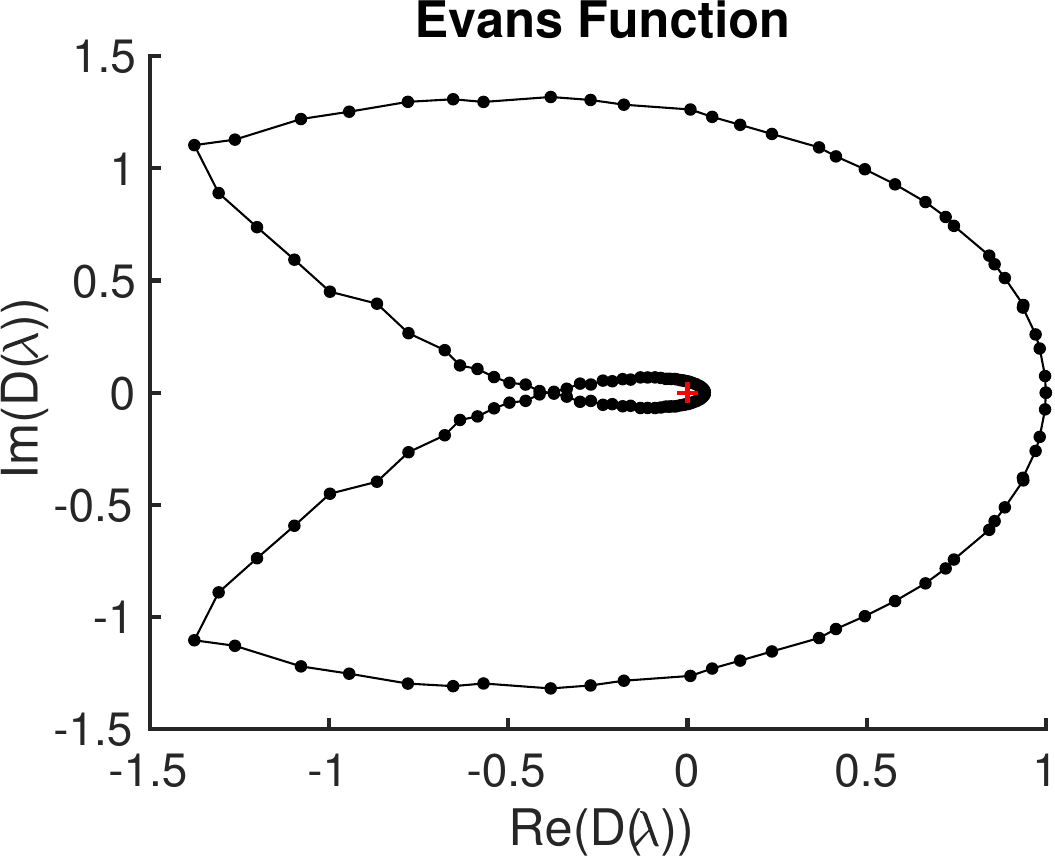} (f)  \includegraphics[scale=0.24]{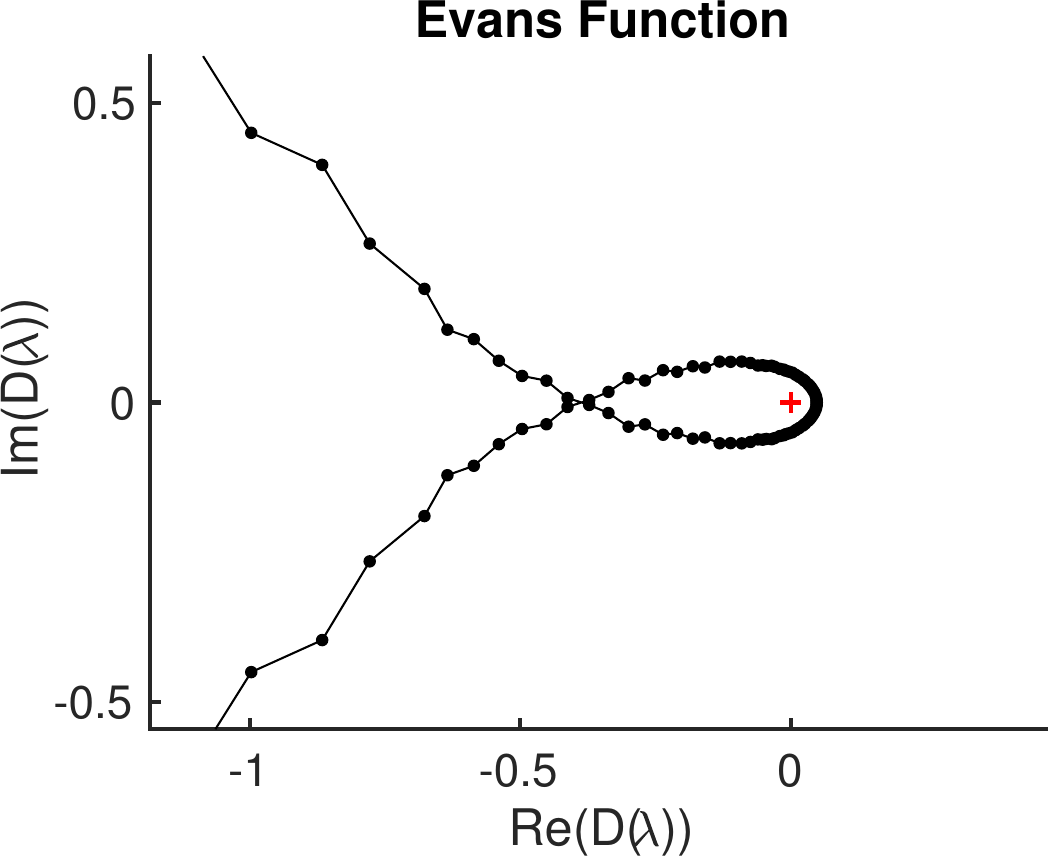}  
\end{array}
\end{equation*}
\end{center}
\caption{The parameters of the local model in this figure are $\mu = 0.5$, $\kappa = 1$, $T_- \approx 1.001$, $T_+ = 2$, $\rho_-  \approx 0.0769$, $\rho_+ = 1$, $u_+ \approx 1.041$, $u_-\approx 13.53$, and $M\approx 1.041$. (a) Plot of $D(\lambda)$ against $\lambda$ where $x_L = -0.5$ and $x_R = 2.15$. (b) Plot of $D(\lambda)$ against $\lambda$ where $x_L = -0.7$ and $x_R = 3.01$. (c) Plot of the whole-line viscous shock profile. (d) Plot of $D(\lambda)$ against $\lambda$ where $x_L = -4.3$ and $x_R = 4.3$. (e) Plot of $\Im(D(\lambda))$ against $\Re(D(\lambda))$  where $x_L = -4.3$, $x_R = 4.3$, and $D(\cdot)$ is evaluated on $\partial ( \{z\in B(0,1e-3):\Re(z)\geq 0\})$. (f) Zoomed-in view of (e).
}
\label{eq:bfz_nonuniqueness}
\end{figure}

 Near the fold curve along which $D(0)=0$, the nullclines of $M_1$ and $M_2$ are by necessity parallel, making 
 computations of their intersection delicate. Our strategy therefore is to move well away 
 from the $ D(0)=0$ curve into negative and positive orientation parts of the plane by varying $(c_1^*,c_2^*)$, then 
 to determine the $M_1$, $M_2$ intersections numerically in this computationally favorable regime 
 where $M_1$ and $M_2$ are substantially transverse.
 We plot the resulting profiles corresponding to $(\hat c_1,\hat c_2)$ and $(\tilde c_1,\tilde c_2)$ in Figures \ref{eq:bfz_multiple2}(a)-(b). We note that there is approximately a 20\% difference between the 
 upper and lower curves over the interval $[-3,3]$ depicted, in terms of the ratio of the 
$\approx 0.2$ maximum difference between the two curves to the $\approx 1.0$ total variation of each curve, 
far more than can be attributed to numerical error. 

Thus, we conclude that these profiles indeed give an explicit example of nonuniqueness.
It may be checked, further, that the two associated parameter pairs as expected lie on opposite sides of the 
fold curve $D(0)\equiv 0$, having value $D(0)$ of opposite signs. 
Indeed, the nullclines of $M_1$ and $M_2$ are approximately parallel, and transverse to the fold curve, 
giving further validation of the numerically observed fold bifurcation scenario.
 
%
%



\subsection{Hopf bifurcation}\label{s:Hopf} Using the same shock parameters in the local model that we used to 
show a bifurcation implying non-uniqueness, but with different choices of left and right boundary,
we can show also the existence of a Hopf-bifurcation. When the finite boundaries are $x_L = -4.3$, $x_R = 4.3$, and the whole-line shock is truncated to $[x_L,x_R]$, the Evans function evaluated on the real line segment $[0,10^{-3}]$ has no zeros, whereas the image of the Evans function evaluated along $\partial ( \{z\in B(0,1e-3):\Re(z)\geq 0\})$ has winding number of two. Thus, there is a complex conjugate pair of eigenvalues with non-zero imaginary part, indicating that a Hopf-bifurcation occurs; see Figure \ref{eq:bfz_nonuniqueness} (d)-(e).

\appendix

\section{Documentation of STABLAB}

In this appendix, we describe additional computational details geared toward the reader interested in reproducing results. In particular, we provide some references regarding the MATLAB-based package STABLAB that we used extensively throughout this paper. STABLAB \cite{Barker_matlab} is a well-tested package for studying stability of traveling waves using the Evans function. This package has been successfully used in a variety of studies; for example see  \cite{Balls-Barker_thesis,BHLynZ1,BHLL,BJNRZ,BLZ,Ghazaryan2020,Lytle,Ozbag2018}.  For an overview of the methods used in STABLAB, please see \cite{Barker_thesis,BHLL}. 

\subsection{Details of feasibility study}\label{appendix:feasible}

To verify the correctness of our code when computing the feasibility set, we independently coded by hand a constant step-size Runge-Kutta-Fehlberg four to fifth order scheme and compared it to the solution we obtained using standard suite software in MATLAB. For improved accuracy, for the large scale study we use MATLAB's \textit{ode15s} \cite{ode15s} routine which is an adaptive step, stiff ODE solver. The solver warnings alert us to finite blowup, and testing the sign of a solution tells us whether or not $u$ and $e$ remain positive throughout the unit interval. In Figure \ref{fig:feasible}, we plot some examples of the feasible set. Note that the feasible set is unbounded (see Lemma \ref{C_unbounded}).

\subsection{Details of the Evans function computations }\label{appendix:evans}

We now provide details about numerical conditioning and algorithm choices for the Evans function computations. For background regarding the methods mentioned, please see \cite{BHLL}. To compute the Evans function, we use the the flux coordinates described in Section 3.1 of \cite{BHLynZ2}, which is equivalent to computing with coordinates $(\rho,u,e,u',e')$ as described in Section \ref{s:evans}. These coordinates are important to use in practice in order to reduce the variation in the image of the Evans function. To improve numerical conditioning of the computation, we evaluate the Evans function wronskian at $x = 1/2$ with ODE solutions given in the definition of the Evans function initialized at $x = 0$ with $\{(0,1,0,0,0)^T,(0,0,1,0,0)^T\}$ and at $x = 1$ with $\{(1,0,0,0,0)^T,(0,1,0,0,0)^T,(0,0,1,0,0)\}$ to recover, as given by Abel's Theorem, a non-vanishing multiple of the Evans function. We also use the method of continuous orthogonalization \cite{HuZ} without the radial equation using Drury's method \cite{Drury} in order to compute the ODE solution, which resolves computational challenges due to differing growth modes.  To verify the correctness of our code, we compute $D(0)$ with the radial equation by initializing the ODE solutions at $x = 0$ only and evolving them to take the determinant at $x = 1$ with the initializing basis there, and check that this matches the value of $D(0)$ computed with the definition given in \eqref{evans}. 

\subsection{Details of the non-uniqueness study}\label{appendix:non_unique}

The boundary value solver we refer to in Section \ref{s:nonuniqueness} is MATLAB's routine bvp5c, which uses a four-stage Lobatto IIIa formula \cite{bvp5c}. We set the tolerance in bvp5c to 1e-6. For the Evans function computations, we use MATLAB's  ode15s with the requested relative and absolute error tolerance set to 1e-10 and 1e-12 respectively. The ode15s routine is a variable-step solver based on variable differentiation formulas of first through fifth orders \cite{ode15s}. 

\subsection{Computational effort}\label{appendix:computational}

Computations were done on a desktop with 128GB Ram and a 4.0GHz i7-6950X Intel processor with 25 MB Cache and 10 cores with 20 threads. Computations were done in Matlab using parallel processing. It took 1.37 days of computation to create the data for the final feasibility study figure, Figure \ref{fig:feasible}. It took  83.8 days of computation on all 10 cores to compute the data for the final Evans function figure, Figure \ref{fig339}. 
Cumulative
computations took longer, exceeding five months. Each of the computations in Section \ref{s:ceg} took a 
substantial part of a day to compute. One of the main reasons the computations were time consuming is stiffness of the 
associated ODE systems. For instance, for the 
ounterexample of Section \ref{s:ceg}, 
continuation of the profile was necessary in order to achieve required accuracy, simple shooting being prohibitively
ill-conditioned.
Indeed, this project is similar numerically in scope and delicacy to those described in \cite{BFZ} and \cite{BMZ}, which together represent a new level of computational challenge in numerical Evans function studies.

\scriptsize
\bibliographystyle{alpha}
\bibliography{biblio}
\normalsize
\end{document}